\documentclass[12pt,twoside,a4paper,leqno,onecolumn]{article}
\usepackage[english,french]{babel}
\usepackage[utf8]{inputenc}
\usepackage[OT1]{fontenc}
\usepackage{sectsty}
\usepackage{fncychap,fancybox}
\usepackage[print,panelleft,gray,paneltoc]{pdfscreen}
\usepackage{mathrsfs,amsthm,amsfonts,amssymb,amscd,amsbsy,amsmath}
\usepackage[retainorgcmds]{IEEEtrantools}
\usepackage{hyperref}
\numberwithin{equation}{section}

\newcommand {\e}{\mathrm{e}}

\newtheorem{theo}{Th\'eor\`eme}

\newtheorem{lem}{Lemme}

\begin{document}
\pagestyle{myheadings}
\markboth{Cherrier, Hanani}{Courbure moyenne prescrite}

\thispagestyle{empty}

\centerline{\begin{Large}\bf Hypersurfaces compactes d'un fibr\'e\end{Large}}

\vskip2mm

\centerline{\begin{Large}\bf vectoriel Riemannien \`a
courbure\end{Large}}

\vskip2mm

\centerline{\begin{Large}\bf moyenne prescrite\end{Large}}

\vskip7mm

\centerline{\textbf{Pascal CHERRIER}\footnote{Current address: Universit\'e de Paris VI, UFR 920 de Math\'ematiques, B.C. 172, 4 place Jussieu, 75252 Paris Cedex 05, France\\ \indent E-mail address : cherrier@ccr.jussieu.fr}}

\vskip2mm

\centerline{et}

\vskip2mm

\centerline{\textbf{Abdellah HANANI}\footnote{Current address: Universit\'e de Lille 1, UFR de Math\'ematiques, B\^at. M2, 59655, Villeneuve d'Ascq Cedex, France\\ \indent E-mail address : abdellah.hanani@math.univ-lille1.fr}}

\vskip8mm

\hrule

\vskip4mm

\noindent\textbf{Abstract. }Let $M$ be a compact Riemannian manifold without boundary and let $E$ be a Riemannian vector bundle over $M$. If $\Sigma $ denotes the sphere subbundle of $E$, we look for embeddings of $\Sigma $ into $E$ admitting a prescribed mean curvature.

\vskip5mm

\noindent\textbf{Mots cl\'es} : Connexions, rel\`evement, courbure moyenne, estim\'ees a priori, les m\'ethodes.

\vskip3mm

\noindent\textbf{Classification mati\`eres 2010} : 35J60, 53C55, 58G30.

\vskip4mm

\hrule


\section{Introduction}

\vskip4mm

Dans cette \'etude, on d\'esigne par $(M,g)$ une vari\'et\'e Riemannienne compacte sans bord, de dimension $n\geq 1$ et $(E,{\tilde g})$ un fibr\'e vectoriel Riemannien sur $M$ de rang $m\geq 2$. On note $\Sigma $ le fibr\'e unitaire correspondant et $E_{*}$ le fibr\'e $E$ priv\'e de la section nulle. On s'int\'eresse \`a la mise en \'evidence d'une hypersurface de $E_{*}$ admettant une courbure moyenne \'egale \`a $K$, une fonction $\mathscr{C}^{\infty }$ strictement positive donn\'ee sur $E_{*}$, et d\'efinie comme la trace de la seconde forme fondamentale relativement \`a la m\'etrique. C'est un probl\`eme global o\`u il s'agit de d\'eterminer un plongement ${\cal Y}$ de $\Sigma $ dans $E_{*}$ ayant une courbure moyenne prescrite.

\vskip4mm

Dans le cadre euclidien, c'est-\`a-dire quand $M$ est r\'eduite \`a un point, un th\'eor\`eme de Bakelman et Kantor [2] assure, en dimension $3$, l'existence d'une telle hypersurface sous la condition que la fonction $K$ d\'ecro\^{\i}t plus vite que la courbure moyenne de sph\`eres concentriques, i.e. il existe deux r\'eels $r_{1}$ et $r_{2}$ tels que $0<r_{1}\leq 1\leq r_{2}$ et
\begin{equation}
K(\xi )>\frac{m-1}{\Vert \xi \Vert }\ \mbox {si}\ \Vert \xi \Vert <r_{1},\
K(\xi )<\frac{m-1}{\Vert \xi \Vert }\ \mbox {si}\ \Vert \xi \Vert >r_{2}
\end{equation}
jointe \`a l'hypoth\`ese de monotonicit\'e
\begin{equation}
\frac{\partial \Big[ rK(r\xi )\Big] }{\partial r}\leq 0,\ \mbox {pour\ tout\ }\xi\in \Sigma .
\end{equation}
Une autre preuve, valable en toute dimension, est donn\'ee par Treibergs et Wei [8] sous les conditions pr\'ec\'edentes. L'hypoth\`ese $(1.2)$ leur a permis d'appliquer la m\'ethode de continuit\'e et leur donne l'unicit\'e \`a homoth\'etie pr\`es.

\vskip4mm

Dans le cadre des fibr\'es envisag\'es ici, des m\'ethodes diff\'erentes de celles des auteurs pr\'ecit\'es s'imposent. En effet, on ram\`ene le probl\`eme \`a la r\'esolution d'une \'equation elliptique sur $\Sigma $ dans laquelle, et contrairement aux cas ci-dessus cit\'es, les d\'eriv\'ees de la fonction inconnue $u$ sont pond\'er\'ees par des puissances de $\e ^u$ dont le degr\'e varie de $0$ \`a $4$. Ceci complique substantiellement l'obtention de l'estim\'ee $\mathscr{C}^1$ et met en p\'eril l'unicit\'e \`a homoth\'etie pr\`es m\^eme sous une hypoth\`ese du type $(1.2)$. 

\vskip4mm

Quand la fonction prescrite ne d\'epend que du rayon, on obtient la caract\'erisation suivante.

\vskip4mm
  
\begin{theo}Soit K une fonction continue sur $E_{*}$, partout strictement
positive, et telle que $K(\xi )=K(\Vert \xi \Vert )$ pour tout $\xi
\in E_{*}$. Alors, il existe un graphe radial $\cal {Y}$ sur $\Sigma
$, de classe $\mathscr{C}^{\infty }$, \`a courbure moyenne \'egale \`a $K$ si et seulement s'il existe un r\'eel $r>0$ tel qu'on ait :
$$K(r\xi )=(m-1)r^{-1}\ quel\ que\ soit\ \xi \in \Sigma .$$
\end{theo}

\vskip5mm

La suffisance dans ce r\'esultat est une cons\'equence de l'expression de l'op\'erateur de Weingarten de l'hypersurface $\Sigma _{r}$, le fibr\'e en sph\`ere de rayon $r>0$. Quant \`a la n\'ec\'essit\'e, elle ne requiert que la continuit\'e de $K$. On g\'en\'eralise ce r\'esultat comme suit.

\vskip4mm

\begin{theo}Soit $K\in \mathscr{C}^{\infty }(E_{*})$ une fonction partout strictement
positive. On fait l'hypoth\`ese qu'il existe deux r\'eels $r_{1}$ et
$r_{2}$, $0<r_{1}\leq 1\leq r_{2}$, tels que les in\'egalit\'es $(1.1)$ soient satisfaites. Notons $\Sigma _{r_{1},r_{2}}=\{\xi \in E;\ r_{1}\leq
\Vert \xi \vert \leq r_{2}\}$ et supposons que $K$ v\'erifie $(1.2)$ en tout point de $\Sigma _{r_{1},r_{2}}$ et que la composante horizontale de son gradient y est partout nulle. Il existe alors un graphe radial $\cal {Y}$ sur $\Sigma $ dont la courbure moyenne est donn\'ee par $K$. Celui-ci \'etant de la forme $\xi \in \Sigma\mapsto e^{u(\xi )}\xi $, o\`u $u\in \mathscr{C}^{\infty }(\Sigma )$ est une fonction \`a gradient horizontal identiquement nul.
\end{theo}
 
\vskip5mm

Dans le but d'\'etablir ce th\'eor\`eme, on constuit le graphe $\cal {Y}$, qui peut \^etre vu comme une fibration dont la fibre ${\cal Y}_{x}$ au dessus du point $x\in M$ est une hypersurface de $E_{x}$ obtenue par projection radiale de la fibre $\Sigma _{x}$ de $\Sigma $, en r\'esolvant une \'equation aux d\'eriv\'ees partielles elliptique sur $\Sigma $. Pour r\'esoudre cette derni\`ere, on utilise l'invariance du degr\'e de
Leray-Schauder d'une homotopie compacte relativement \`a un domaine born\'e du sous-ensemble des fonctions de classe $\mathscr{C}^{1,\alpha }(\Sigma )$ dont la composante horizontale du gradient est identiquement nulle. L'application de celle-ci repose en partie sur l'obtention d'une estim\'ee a priori dans $\mathscr{C}^{1,\alpha }(\Sigma )$. Remarquons que, dans ce r\'esultat, l'hypoth\`ese $(1.2)$ est faite pour assurer
la nullit\'e du gradient horizontal de la solution. Celle-ci peut \^etre omise. A cet \'egard, on d\'emontre le th\'eor\`eme suivant.

\vskip4mm

\begin{theo}Soit $K\in \mathscr{C}^{\infty }(E_{*})$ une fonction partout strictement positive. On suppose qu'il existe deux r\'eels $r_{1}$ et $r_{2}$, $0<r_{1}\leq 1\leq r_{2}$, tels que les in\'egalit\'es $(1.1)$ soient satisfaites. Il existe alors un graphe radial ${\cal Y}$ sur $\Sigma $, de classe $\mathscr{C}^{\infty }$, dont la
courbure moyenne est donn\'ee par $K$, et tel que $r_{1}\leq \Vert \xi \Vert
\leq r_{2}$ pour tout $\xi \in {\cal Y}$.
\end{theo}

\vskip5mm

Pour prouver l'existence d'une solution sous ces conditions, on se place dans le cadre fonctionnel $\mathscr{C}^{\infty }$ et on applique alors le th\'eor\`eme de point fixe de Nagumo [7]. Ceci est due au fait qu'en vue d'\'etablir l'estim\'ee $\mathscr{C}^{1}$ n\'ec\'essaire pour r\'esoudre, il faut disposer auparavant d'un contr\^ole sur le laplacien horizontal de toute solution \'eventuelle. Dans ce but, on \'etend l'\'equation \`a r\'esoudre au voisinage de $\Sigma $, les valeurs \`a prescrire au bord \'etant du type Neumann. L'\'equation \`a consid\'erer au voisinage de $\Sigma $ est de sorte qu'une d\'erivation radiale et restriction \`a $\Sigma $ de celle-ci donne le contr\^ole souhait\'e.

\vskip4mm

Nous pr\'esentons cet article en quatre parties. A la derni\`ere, on donne la preuve des th\'eor\`emes et une remarque montrant que l'hypoth\`ese de croissance des th\'eor\`emes 2 et 3 est dans un certain sens la meilleure possible. L'estim\'ee a priori de la norme $\mathscr{C}^{1}$ n\'ec\'essaire pour r\'esoudre est pr\'esent\'ee \`a la troisi\`eme partie. Une mise en \'equation est donn\'ee \`a la seconde partie et, pour plus de clart\'e, on consacre la premi\`ere partie de cet article \`a des notations et rappels pr\'eliminaires. 

\vskip6mm

\section{Notations et pr\'eliminaires}

\vskip4mm
 
\noindent \textbf{1- }Soit $(M,g)$ une vari\'et\'e Riemannienne compacte sans bord de dimension $n\geq 1$. Soient $(E,{\tilde g})$ un fibr\'e vectoriel Riemannien sur $M$ de rang $m\geq 2$, $\pi $ la projection naturelle de $E$ sur $M$ et $E_{*}$ le fibr\'e $E$ priv\'e de la section nulle. On note $\nabla $ la connexion de Levi-Civita de la vari\'et\'e $(M,g)$ et ${\tilde \nabla }$ une connexion m\'etrique sur le fibr\'e $(E,{\tilde g})$.

\vskip4mm

Soit $U$ un ouvert de $M$, domaine du syst\`eme de coordonn\'ees $(x^{i})_{1\leq
i\leq n}$ sur $M$, et au dessus duquel $E$ est trivial. Notons $\displaystyle \varepsilon _{i}=\frac{\partial }{\partial x^i}$, $1\leq i\leq n$, et consid\'erons un rep\`ere de sections de $E$ au dessus de $U$ not\'e $\displaystyle (\varepsilon _{\alpha })$, $n+1\leq \alpha \leq n+m$. Si $\xi \in \pi ^{-1}(U)$ et $x=\pi (\xi )$, on \'ecrit $\xi =y^{\alpha }\varepsilon _{\alpha }(x)$; $(x^{i},y^{\alpha })_{i,\alpha }$ est alors un syst\`eme de coordonn\'ees sur $\pi ^{-1}(U)$.

\vskip4mm

Soient $x\in U$, $\displaystyle \xi =\xi ^{\alpha }\varepsilon _{\alpha }\in \pi
^{-1}(x)$ et $\displaystyle X =X^{i}\varepsilon _{i}$ un champ de vecteurs tangents \`a
$U$. Le rel\`evement horizontal de $X_{\vert _{x}}$ en $\xi $ est l'\'el\'ement
$X^{H}(\xi )\in T_{\xi }E$ d\'efini par
\begin{equation}
X^{H}(\xi )=X^{i}(x)\left(\frac{\partial }{\partial x^i}_{\vert _{\xi }}-\xi ^{\alpha }\Gamma ^{\beta }_{i\alpha }\frac{\partial }{\partial y^{\beta }}_{\vert_{\xi }}\right),\end{equation}
o\`u les $\Gamma ^{\beta }_{i\alpha }$, $i\in \{1,...,n\}$ et $\alpha ,\beta
\in \{n+1,...,n+m\}$ d\'esigent les symboles de Christoffel de la connexion
${\tilde \nabla }$ dans le syst\`eme de coordonn\'ees $(x^{i},y^{\alpha })$.

\vskip4mm

Sur l'ouvert $\pi ^{-1}(U)$, on consid\`ere alors le rep\`ere mobile  
$${\cal S}=\{e_{i},e_{\alpha }/\ i=1,...,n\ \mbox {et}\ \alpha =n+1,...,n+m\},$$
o\`u $\displaystyle e_{\alpha }=\frac{\partial }{\partial y^{\alpha }}$ et
$e_{i}$ est le rel\`evement horizontal de $\displaystyle \varepsilon _{i}$. Ainsi
\begin{equation}
e_{i}=\frac{\partial }{\partial x^i}-y^{\alpha }\Gamma ^{\beta }_{i\alpha}\frac{\partial }{\partial y^{\beta }}.
\end{equation}

\vskip2mm

On d\'efinit sur la vari\'et\'e $E$ une m\'etrique Riemannienne $G$ en posant
\begin{equation}
G(e_{i},e_{j})=g(\epsilon _{i},\epsilon _{j}),\ \ G(e_{\alpha },e_{\beta
})={\tilde g}(\epsilon _{\alpha },\epsilon _{\beta }),\ \ G(e_{i},e_{\alpha })=0.
\end{equation}
On dispose alors, sur la vari\'et\'e $E$, de la connexion de Levi-Civita associ\'ee \`a la m\'etrique $G$ et de la connexion $D$ de Sasaki [9] et [10] d\'efinie par
\begin{equation}
D_{e_{i}}e_{j}=\Gamma ^{k}_{ij}e_{k},\ D_{e_{i}}e_{\alpha }=\Gamma ^{\beta
}_{i\alpha }e_{\beta },\ D_{e_{\alpha }}e_{i}=D_{e_{\alpha }}e_{\beta }=0,
\end{equation} 
o\`u les $\Gamma ^{k}_{ij}$, $i,j,k\in \{1,...,n\}$, d\'esignent les symboles de
Christoffel de la connexion $\nabla $ dans le syst\`eme de coordonn\'ees $(x^{i})$. La connexion $D$ est $G$-m\'etrique et ne co\"{\i}ncide pas avec la connexion de Levi-Civita de $G$; la torsion $T$ de $D$ est non nulle. En effet, notant $\mathsf{S}$
le tenseur de courbure de ${\tilde \nabla }$, les composantes de $T$ dans le rep\`ere mobile ${\cal S}$ sont toutes nulles sauf
\begin{equation}
T^{\alpha }_{ij}=-y^{\beta }S^{\alpha }_{\beta ij}.
\end{equation}
Les composantes dans ${\cal S}$ du tenseur de courbure ${\cal R}$ de $D$ sont donn\'ees par
$$R_{dcab}=G\left((D_{e_{a}e_{b}}-D_{e_{b}e_{a}}-D_{[e_{a},e_{b}]})e_{c},e_{d}\right)$$
et 
$$R^{d}_{cab}=G^{de}R_{ecab}.$$
Par un calcul direct, on montre qu'en particulier, pour $1\leq i,j\leq n$ et
$n+1\leq \alpha ,\beta ,\lambda ,\mu \leq n+m$, on a :
\begin{equation}
R^{i}_{\alpha \ \beta \ j}=R^{\lambda }_{\alpha \ \beta \ j}=R^{i}_{\alpha \ \beta \ \mu }=R^{\lambda }_{\alpha \ \beta \ \mu }\equiv 0.
\end{equation}

\vskip5mm

\noindent \textbf{2- }On note $\Sigma =\{\xi \in E/\ \Vert \xi \Vert =1\}$ le sous-fibr\'e unitaire de $E$, $\pi_{1} $ la surjection naturelle du fibr\'e $\Sigma $ sur $M$, $r$ la fonction $r(\xi )=\Vert \xi \Vert $ et $\nu $ le champ radial unitaire. Sur l'ouvert $\pi ^{-1}(U)$ muni des coordonn\'ees $(x^{i},y^{\alpha })$, $1\leq i\leq n$ et $n+1\leq \alpha \leq n+m$, le champ $\nu $ est donn\'e par
$$\nu =r^{-1}y^{\alpha }\frac{\partial }{\partial y^{\alpha }}$$ 
et l'ouvert $\pi _{1}^{-1}(U)$ de $\Sigma $ est donn\'e, comme sous-vari\'et\'e de dimension $n+m-1$ de $\pi ^{-1}(U)$, par
$${\tilde g}_{\alpha \beta }y^{\alpha }y^{\beta }=1.$$
Ce qui, par un calcul direct montre que le champ radial unitaire $\nu $ est normal \`a $\Sigma $. De l'expression de $\nu $, celle du rel\`evement horizontal d'un vecteur $X\in T_{x}M$ au point $\xi =(x,y^{\alpha })\in \Sigma $, donn\'ee par $(2.2)$, et de la d\'efinition de la m\'etrique $G$, on voit que l'espace tangent \`a $\Sigma $ au point $\xi $ est une somme directe du sous-espace horizontal $H_{\xi }E$ de $T_{\xi}E$ et de l'espace tangent \`a la fibre de $\Sigma $ passant par $\xi $.

\vskip4mm

Dans tout ce qui suit, le param\`etre $\mu _{a}$ sera \'egal \`a $0$ ou $1$ selon que la direction $a$ est verticale ou horizontale. Fixons un rep\`ere mobile tangent \`a $E$ de la forme
$${\cal {R}}=\{e_{i},e_{\alpha }/\ i=1,...,n\ \mbox{et}\ \alpha
=n+1,...,n+m\},$$
o\`u les $e_{i}$ sont des champs de vecteurs horizontaux obtenus par rel\`evement horizontal d'un rep\`ere mobile $(\epsilon _{i})_{1\leq i\leq n}$ sur $M$ et o\`u les $e_{\alpha }$ sont des champs de vecteurs verticaux avec $e_{n+m}=\nu $. On notera
$$\displaystyle {\cal {R}}^{*}=\{\omega ^{A},\ A\leq n+m \}$$
le corep\`ere dual de ${\cal {R}}$. La restriction \`a $\Sigma $ des champs $e_{a}$, pour $a\leq n+m-1$, constitue un rep\`ere mobile ${\cal {R}}_{1}$ tangent \`a $\Sigma $ et
$${\cal {R}}_{1}^{*}=\{\sigma ^{a},\ a\leq n+m-1 \},$$
o\`u $\sigma ^{a}$ est l'image inverse de $\omega ^{a}$ par $i$, injection canonique de $\Sigma $ dans $E$, n'est autre que le corep\`ere dual de ${\cal {R}}_{1}$. 

\vskip4mm

Appliquons $D$ \`a $e_{a}$, le r\'esultat est une 1-forme sur $E$ \`a valeurs dans $TE$ dont l'expression dans la rep\`ere ${\cal {R}}$ nous permet d'introduire la matrice $(\omega ^{A}_{B})$ de 1-formes d\'efinie par les \'egalit\'es
\begin{equation}
De_{A}=\omega ^{B}_{A}\otimes e_{B}.
\end{equation}
Du fait que $\nu $ est unitaire et puisque $D$ est $G$-m\'etrique, on voit que
$$\displaystyle G(D_{e{_{a}}}\nu ,\nu )=0,\ \mbox{pour}\ 1\leq a\leq n+m.$$
On en d\'eduit que, sur $\Sigma _{r}$, le fibr\'e en sph\`eres de rayon $r$,
\begin{equation}
\omega ^{n+m}_{n+m}=0.
\end{equation}
D'autre part, si $X$ est un champ de vecteurs verticaux tangents \`a $\Sigma _{r}$ et donc v\'erifiant $\displaystyle i_{*}Xr=0$ sur $\Sigma _{r}$, on aura
\begin{equation}
D_{i_{*}X}\nu =r^{-1}(i_{*}X.y^{\alpha })\frac{\partial }{\partial y^{\alpha }}+r^{-1}y^{\alpha }D_{i_{*}X}\left(\frac{\partial }{\partial y^{\alpha }}\right)=r^{-1}i_{*}X.
\end{equation}
Eu \'egard \`a $(2.4$) qui implique que $\displaystyle D_{i_{*}X}\left(\frac{\partial  }{\partial y^{\alpha }}\right)=0$. De m\^eme, compte tenu de la d\'efinition $(2.4)$ de la connexion $D$ et de l'expression $(2.2)$ du rel\`evement horizontal d'un champ de vecteurs, si dans $U$ , $\displaystyle \varepsilon _{i}=\varepsilon _{i}^{j}\frac{\partial }{\partial x^j}$ et $e_{i}=(\varepsilon _{i})^{H}$, on voit que
$$D_{e_{j}}\nu =r^{-1}\varepsilon _{j}^{h}\left(\frac{\partial y^{\alpha }}{\partial x^h}-y^{\lambda }\Gamma ^{\mu }_{h\lambda }\frac{\partial y^{\alpha }}{\partial y^{\mu }}\right)\frac{\partial }{\partial y^{\alpha }}+r^{-1}y^{\alpha }\varepsilon ^{h}_{j}\Gamma ^{\mu }_{h\alpha }\frac{\partial }{\partial y^{\mu }}.$$
D'o\`u
\begin{equation}
D_{e_{j}}\nu =r^{-1}\epsilon _{j}^{h}\left[-y^{\lambda }\Gamma ^{\mu }_{h\lambda }\delta ^{\alpha }_{\mu }\frac{\partial }{\partial y^{\alpha }}+y^{\alpha }\Gamma ^{\mu }_{h\alpha }\frac{\partial }{\partial y^{\mu }}\right]=0.
\end{equation}
Ainsi, compte tenu de $(2.9)$ et $(2.10)$,
\begin{equation}
D_{e_{a}}\nu =(1-\mu _{a})r^{-1}e_{a}\quad \mbox{pour}\quad a\leq n+m-1.
\end{equation}
Reportons dans $(2.7)$, il en d\'ecoule que, sur $\Sigma _{r}$, on a 
\begin{equation}
\omega ^{a}_{n+m}=(1-\mu _{a})r^{-1}\omega ^{a}\ \mbox{pour}\ a\leq n+m-1
\end{equation}
et par suite, compte tenu de $(2.8)$,
\begin{equation}
D_{\nu }\nu =0.
\end{equation}
D'autre part, la compatibilit\'e de la connexion $D$ avec la m\'etrique $G$ montre que
$$G(D_{e_{b}}e_{a},\nu )=-G(e_{a},D_{e_{b}}\nu )$$
ce qui, compte tenu de $(2.12)$, implique que
\begin{equation}
\omega ^{n+m}_{a}(e_{b})=-(1-\mu _{b})r^{-1}G_{ab}\ {\rm pour\it \ }a,b\leq n+m-1.
\end{equation}

\vskip2mm

A pr\'esent, on donne les composantes du tenseur de courbure ${\tilde {\cal R}}$ de $\Sigma $. Utilisons l'\'equation de Gauss et $(2.14$) on voit que
$${\tilde R}_{dcab}=R_{dcab}+(1-\mu _{a})(1-\mu _{b})(G_{ad}G_{bc}-G_{ac}G_{bd}).$$
Ainsi, tenons compte de $(2.6)$, nous obtenons la valeur des composantes qui seront utilis\'ees ult\'erieurement:
\begin{equation}
{\tilde R}^{j}_{\alpha \beta \gamma }={\tilde R}^{j}_{\alpha \beta i}={\tilde R}^{\gamma }_{\alpha \beta i}=0,\ n+1\leq \alpha ,\beta ,\gamma \leq n+m-1\ \mbox{et}\ 1\leq i,j\leq n,
\end{equation} 
et
\begin{equation}
{\tilde R}^{\lambda }_{\alpha \beta \mu }=\delta ^{\lambda }_{\beta
}G_{\alpha \mu }-\delta ^{\lambda }_{\mu }G_{\alpha \beta },\ n+1\leq \alpha
,\beta ,\lambda ,\mu \leq n+m-1.
\end{equation}

\vskip8mm

\section{Mise en \'equation}

\vskip4mm

Dans cette partie, on donnera l'\'equation de la courbure moyenne. On conserve toutes les notations du dernier paragraphe de la section pr\'ec\'edente. Soit $u\in \mathscr{C}^{2}(\Sigma )$ une fonction qu'on prolonge \`a $E_{*}$ en la maintenant radialement constante. On consid\`ere l'application ${\cal Y}$ de $\Sigma $ dans $E$ telle que 
$${\cal Y}(\xi )=e^{u(\xi )}\xi ,\ \ {\rm pour\ \it }\xi \in \Sigma .$$
La diff\'erentielle $D{\cal Y}$ de ${\cal Y}$ est une section de $T^{*}\Sigma \otimes TE$ qu'on peut \'ecrire sous la forme 
$$D{\cal Y}=\sigma ^{a}\otimes E_{a},\ {\rm o\grave u\ pour\it \ }a\leq n+m-1,\ E_{ a}=D{\cal Y}(e_{a}).$$
Les champs de vecteurs $\{E_{a},a\leq n+m-1\}$ d\'efinissent un rep\`ere mobile tangent au graphe ${\cal Y}$. Un calcul direct donne 
$$E_{a}=e_{a}+e^{u}D_{a}u\nu $$ 
et le champ unitaire normal \`a ${\cal Y}$ est donc
$${\tilde \nu }=f(\nu -e^{u}D^{a}ue_{a}),\ \ f=(1+e^{2u}D_{a}uD^{a}u)^{-{\frac{1}{2}}}.$$
La m\'etrique $H$ induite par $G$ sur ${\cal Y}$ est
$$H=H_{ab}\sigma ^{a}\otimes \sigma ^{b},\ {\rm o\grave u\it \ }H_{ab}=G(E_{a},E_{b})=G_{ab}+e^{2u}D_{a}uD_{b}u$$
et la r\'esolution de l'\'equation $H^{ab}H_{bc}=\delta ^{a}_{c}$ donne les composantes
contravariantes de celle-ci, on v\'erifie que : 
$$H^{ab}=G^{ab}-f^{2}e^{2u}D^{a}uD^{b}u.$$
Enfin, rappelons qu'avec le choix ci-dessus du champ normal unitaire ${\tilde \nu }$, la courbure moyenne ${\cal M}_{_{\cal Y}}$ de ${\cal Y}$ est donn\'ee par
$${\cal M}_{_{\cal Y}}=H^{ab}G(D_{E_{a}}{\tilde \nu },E_{b}).$$
Dans le corep\`ere ${\cal R}^*$, la diff\'erentielle de la fonction $u$ est donn\'e par 
$$\displaystyle du=\sum _{a=1}^{n+m-1}D_{a}u\omega ^{a}.$$
La composante $D_{a}u$ est homog\`ene de degr\'e $0$ ou $-1$ selon que la direction $a$ est horizontale ou verticale, i.e. $D_{a}u$ est homog\`ene de degr\'e $(\mu _{a}-1)$. De m\^eme, on a
$$D_{ab}u=D^2u(e_{a},e_{b})=(D_{e_{a}}Du)(e_{b}).$$
D'o\`u
$$D_{ab}u=D_{e_{a}}\Big(Du(e_{b})\Big)-Du(D_{e_{a}}e_{b})$$
et on v\'erifie que la composante $\displaystyle D_{ab}u$ est homog\`ene de degr\'e $(\mu _{a}+\mu _{b}-2)$. En particulier, sur ${\cal Y}$, on peut \'ecrire 
$$D_{a\nu }u=D_{e_{a}}\Big(Du(\nu )\Big)-Du(D_{e_{a}}\nu )=-Du(D_{e_{a}}\nu ).$$
La derni\`ere \'egalit\'e d\'ecoule du fait que $u$ est une fonction radialement constante. La relation $(2.11)$ implique alors que
\begin{equation}
D_{a\nu }u=-(1-\mu _{a})e^{-u}D_{a}u\ {\rm pour\it \ }a\leq n+m-1.
\end{equation}
Tenons compte de la relation $(2.13)$, un calcul analogue montre que
\begin{equation}
D_{\nu \nu }u=0.
\end{equation}
A pr\'esent, de la nature du prolongement de $u$ qui implique que $D_{\nu }u=0$, et la
relation $(2.13)$ qui dit que $\displaystyle D_{\nu }\nu =0$, on voit que
$$D_{E_{a}}E_{b}=D_{e_{a}}(e_{b}+e^uD_{b}u\nu )+e^uD_{a}D_{\nu }e_{b}+e^{2u}D_{a}u(\nu .D_{b}u)\nu .$$
Le fait que les composantes de la torsion de la forme $T^{a}_{n+m\ b}$ sont identiquement nulles et la relation $(2.11)$ permettent d'\'ecrire :  
$$D_{\nu }e_{b}=D_{e_{b}}\nu =(1-\mu {b})e^{-u}e_{b}.$$ 
D'autre part,
$$D_{\nu }(D_{b}u)=D_{\nu b}u+Du(D_{\nu }e_{b})=-(1-\mu _{b})e^{-u}D_{b}u+Du\Big((1-\mu
_{b})e^{-u}e_{b}\Big).$$
La derni\`ere \'egalit\'e est une cons\'equence de $(3.1)$ et la relation ci-dessus. Ainsi, on voit que $D_{\nu }(D_{b}u)=0$ et donc
$$D_{E_{a}}E_{b}=D_{e_{a}}e_{b}+e^{u}\Big[D_{e_{a}}(D_{b}u)+D_{a}uD_{b}u\Big]\nu
+e^{u}D_{b}uD_{e_{a}}\nu +(1-\mu _{b})D_{a}ue_{b}.$$ 
Une \'egalit\'e qui, compte tenu de $(2.11)$, s'\'ecrit sous la forme
\begin{equation}
D_{E_{a}}E_{b}=D_{e_{a}}e_{b}+e^{u}\Big[D_{e_{a}}(D_{b}u)+D_{a}uD_{b}u\Big]\nu
+D_{b}u(1-\mu _{a})e_{a}+D_{a}u(1-\mu _{b})e_{b}.
\end{equation}
Or,
$$D_{e_{a}}e_{b}=\sum _{1}^{n+m-1}\omega ^{c}_{b}(e_{a})e_{c}+\omega ^{n+m}_{b}(e_{a})\nu =\sum _{1}^{n+m-1}\omega ^{c}_{b}(e_{a})e_{c}-(1-\mu _{a})e^{-u}G_{ab}\nu $$
Donc, tenons compte de $(2.14)$, on obtient
$$D_{e_{a}}e_{b}=\sum _{1}^{n+m-1}\omega ^{c}_{b}(e_{a})e_{c}-(1-\mu
_{a})e^{-u}G_{ab}\nu .$$
Par suite, la d\'efinition de la d\'eriv\'ee covariante permet d'\'ecrire 
$$D_{e_{a}}(D_{b}u)=D_{ab}u+Du(D_{e_{a}}e_{b})=D_{ab}u+\sum _{1}^{n+m-1}\omega
^{c}_{b}(e_{a})D_{c}u.$$
De sorte que
$$D_{e_{a}}e_{b}+e^{u}(e_{a}D_{b}u)\nu =\sum _{1}^{n+m-1}\omega ^{c}_{b}(e_{a})E_{c}+e^{-u}\Big[e^{2u}D_{ab}u-(1-\mu _{a})G_{ab}\Big]\nu.$$
Reportons dans (3.3), compte tenu du fait que $\displaystyle E_{a}=e_{a}+e^{u}D_{a}u\nu $, la relation qui en r\'esulte peut s'\'ecrire sous la forme :
$$\begin{array}{ccl}\displaystyle D_{E_{a}}E_{b}&=&\displaystyle \sum _{1}^{n+m-1}\omega ^{c}_{b}(e_{a})E_{c}+D_{b}u(1-\mu _{a})E_{a}+D_{a}u(1-\mu _{b})E_{b}\\ &&\displaystyle +e^{-u}\Big[-(1-\mu _{a})H_{ab}+\mu _{b}e^{2u}D_{a}uD_{b}u+e^{2u}D_{ab}u\Big]\nu \end{array}$$
et par suite, eu \'egard au fait que $G(E_{c},{\tilde \nu })=0$ pour $c\leq n+m-1$, on obtient : 
$$G(D_{E_{a}}E_{b},{\tilde \nu })=fe^{-u}\Big[-(1-\mu _{a})H_{ab}+\mu _{b}e^{2u}D_{a}uD_{b}u+e^{2u}D_{ab}u\Big].$$
Ainsi et puisque $G(D_{E_{a}}{\tilde \nu },E_{b})=-G(D_{E_{a}}E_{b},{\tilde \nu })$, on trouve que 
$$f^{-1}e^{u}{\cal M}_{_{\cal Y}}=H^{ab}\Big[(1-\mu _{a})H_{ab}-\mu
_{b}e^{2u}D_{a}uD_{b}u-e^{2u}D_{ab}u\Big].$$
Prenons l'image inverse sur $\Sigma $ de cette \'equation, tenons compte de l'homogen\'eit\'e des d\'eriv\'ees covariantes de $u$ et notons
${\tilde D}^{a}u=e^{\mu _{a}u}D^{a}u$ et ${\tilde D}^{ab}u=e^{(\mu _{a}+\mu _{b})u}D^{ab}u$, on obtient
$$A^{ab}(u){\tilde D}_{ab}u=-v_{2}+(m-1)(1+v)-(1+v)^{\frac{3}{2}}e^{u}{\cal M}_{_{\cal Y}}(e^{u}\xi ),$$
o\`u $\displaystyle A^{ab}(u)=(1+v)G^{ab}-{\tilde D}^{a}u{\tilde D}^{b}u$ et
$$v=v_{1}+v_{2},\ v_{1}=(1-\mu _{a}){\tilde D}_{a}u{\tilde D}^{a}u,\ v_{2}=\mu
_{a}{\tilde D}_{a}u{\tilde D}^{a}u.$$
Ainsi, la recherche d'un graphe radial \`a courbure moyenne prescrite $K$ revient \`a r\'esoudre sur $\Sigma $ l'\'equation elliptique suivante :
\begin{equation}
A^{ab}(u){\tilde D}_{ab}u=-v_{2}+(m-1)(1+v)-(1+v)^{\frac{3}{2}}e^{u}K(e^{u}\xi ).
\end{equation}

\vskip6mm

\section{Estimations a priori}

\vskip4mm

\begin{lem}Soient $K\in \mathscr{C}^{\infty }(E_{*})$ une fonction partout
strictement positive et $u\in \mathscr{C}^{3}(\Sigma )$ une solution de l'\'equation
\begin{equation}
A^{ab}{\tilde D}_{ab}u=-v_{2}+(m-1)(1+v)-(1+v)^{\frac{3}{2}}e^{u}K(e^{u}\xi
)\equiv F
\end{equation}  
telle qu'il existe deux r\'eels strictement positifs $r_{1}$ et $r_{2}$ v\'erifiant
$r_{1}\leq e^{u}\leq r_{2}$. Notons $\Sigma _{r_{1},r_{2}}=\{\xi \in E;\ r_{1}\leq \Vert \xi \Vert \leq r_{2}\}$ et $r$ la fonction $r(\xi )=\Vert \xi \Vert $, et supposons qu'en tout point de $\Sigma $ o\`u $v\geq 1$, on ait
\begin{equation}
(\mu _{a}+\mu _{b})G^{ab}{\tilde D}_{ab}u\leq -\alpha \sqrt{v},
\end{equation}
o\`u l'on a not\'e
\begin{equation}
\alpha =2\sup _{\xi \in \Sigma }\left[\sup _{\Sigma _{r_{1},r_{2}}}\left\vert r\frac{\partial[rK(r\xi )]}{\partial r}\right\vert\right].
\end{equation}
Alors il existe une constante positive $C_{1}$ ne d\'ependant que de $(M,g)$, $(E,{\tilde g})$, $r_{1}$, $r_{2}$, $\displaystyle \max _{\Sigma _{r_{1},r_{2}}}K$ et $\Vert K\Vert _{\mathscr{C}^{1}(\Sigma _{r_{1},r_{2}})}$ telle que $\vert Du\vert \leq C_{1}$ partout dans $\Sigma $.
\end{lem}

\vskip5mm

\noindent\textit{D\'emonstration. }Soient $u\in \mathscr{C}^{3}(\Sigma )$ une solution de $(4.1)$, $l$ un r\'eel strictement positif fix\'e ult\'erieurement et $\Gamma $ la fonctionnelle d\'efinie, sur $\Sigma $, par
\begin{equation}
\Gamma (u)=(1+v)e^{-lu}.
\end{equation}
En un point $\xi \in \Sigma $ o\`u $\Gamma (u) $ atteint son maximum, on a :
\begin{equation}
\frac{D_{a}\Gamma }{\Gamma }=\frac{D_{a}v}{1+v}-lD_{a}u=0
\end{equation}
et
$$A^{ab}{\tilde D}_{ab}\Gamma \leq 0$$
c'est-\`a-dire, tenons compte de $(4.1)$,
\begin{equation}
A^{ab}{\tilde D}_{ab}v-\frac{A^{ab}{\tilde D}_{a}v{\tilde D}_{b}v}{1+v}-l(1+v)F\leq
0.
\end{equation}
Or,
$$\begin{array}{c}\displaystyle A^{ab}{\tilde D}_{ab}v=4v_{2}A^{ab}{\tilde D}_{a}u{\tilde D}_{b}u+2v_{2}A^{ab}{\tilde D}_{ab}u+8A^{ab}\mu _{c}{\tilde D}_{ac}u{\tilde D}^{c}u{\tilde D}_{b}u\\ \\ \displaystyle +2A^{ab}{\tilde D}_{a}^{c}u{\tilde D}_{bc}u+2A^{ab}{\tilde D}_{abc}u{\tilde D}^{c}u,\end{array}$$
o\`u l'on a not\'e ${\tilde D}_{abc}u=e^{(\mu _{a}+\mu _{b}+\mu _{c})u}D_{abc}u$. Une permutation de l'ordre des indices de d\'erivation covariante dans le terme des d\'eriv\'ees troisi\`emes, celle-ci g\`en\`ere des termes en torsion et en courbure, montre que
\begin{equation}
\begin{array}{c}A^{ab}{\tilde D}_{ab}v=4v_{2}v+2v_{2}F+8A^{ab}\mu _{c}{\tilde D}_{ac}u{\tilde D}^{c}u{\tilde D}_{b}u\\ \\ \displaystyle
+2A^{ab}{\tilde D}_{a}^{c}u{\tilde D}_{bc}u+2A^{ab}{\tilde D}_{cab}u{\tilde D}^{c}u+4E_{1}+E_{2},\end{array}
\end{equation}
o\`u les termes $E_{1}$ et $E_{2}$ sont donn\'es par
$$E_{1}=-A^{ab}e^{(\mu _{a}+\mu _{b}+\mu _{c})u}T^{d}_{ac}D_{bd}u{\tilde D}^{c}u$$
et
$$E_{2}=-2A^{ab}e^{(\mu _{a}+\mu _{b}+\mu _{c})u}\Big({\tilde R}^{d}\
_{bac}+D_{a}T^{d}_{bc}-T^{h}_{ac}T^{d}_{hb}\Big) D_{d}u{\tilde D}^{c}u.$$

\vskip2mm

A pr\'esent, on d\'erive une fois l'\'equation $(4.1)$, il vient
$$A^{ab}{\tilde D}_{cab}u+(\mu _{a}+\mu _{b})A^{ab}{\tilde D}_{ab}u{\tilde D}_{c}u+({\tilde D}_{c}A^{ab}){\tilde D}_{ab}u={\tilde D}_{c}F.$$
D\'eveloppons $({\tilde D}_{c}A^{ab})$ et saturons par ${\tilde D}^{c}u$, on obtient :
\begin{equation}
\begin{array}{c}A^{ab}{\tilde D}_{cab}u{\tilde D}^{c}u={\tilde D}_{c}F{\tilde D}^{c}u-v(\mu _{a}+\mu _{b})A^{ab}{\tilde D}_{ab}u-{\tilde D}_{c}v{\tilde D}^{c}uG^{ab}{\tilde D}_{ab}u\\ \\ \displaystyle +v(\mu _{a}+\mu _{b}){\tilde D}^{a}u{\tilde D}^{b}u{\tilde D}_{ab}u+{\tilde D}_{c}\ ^{a}u{\tilde D}^{c}u{\tilde D}^{b}u{\tilde D}_{ab}u+{\tilde D}^{a}u{\tilde D}^{c}u{\tilde D}_{c}\ ^{b}u{\tilde D}_{ab}u.\end{array}
\end{equation}
D'apr\`es $(4.1)$, on peut \'ecrire $\displaystyle G^{ab}{\tilde D}_{ab}u=\frac{F+{\tilde D}^{a}u{\tilde D}^{b}u{\tilde D}_{ab}u}{1+v}$. Il s'en suit que
$${\tilde D}_{c}v{\tilde D}^{c}u{\tilde D}_{a}^{a}u=\frac{F}{1+v}{\tilde D}_{c}v{\tilde D}^{c}u+\frac{{\tilde D}^{a}u{\tilde D}^{b}u{\tilde D}_{ab}u{\tilde D}_{c}v{\tilde D}^{c}u}{1+v}.$$
D\'eveloppons ${\tilde D}_{c}v$ dans le second terme du membre de droite, il vient
$${\tilde D}_{c}v{\tilde D}^{c}u{\tilde D}_{a}^{a}u=\frac{F{\tilde D}_{c}v{\tilde D}^{c}u}{1+v}+\frac{2\Big( {\tilde D}^{a}u{\tilde D}^{b}u{\tilde D}_{ab}u\Big)^{2}}{1+v}+\frac{2v_{2}v}{1+v}{\tilde D}^{a}u{\tilde D}^{b}u{\tilde D}_{ab}u.$$
Or
\begin{equation}
{\tilde D}_{a}v=2{\tilde D}^{b}u{\tilde D}_{ab}u+2v_{2}{\tilde D}_{a}u.
\end{equation}
Il en d\'ecoule que
\begin{equation}
\begin{array}{c}\displaystyle {\tilde D}_{c}v{\tilde D}^{c}uG^{ab}{\tilde D}_{ab}u=\frac{F}{1+v}{\tilde D}_{c}v{\tilde D}^{c}u+\frac{2}{1+v}\Big( {\tilde D}^{a}u{\tilde D}^{b}u{\tilde D}_{ab}u\Big)^{2}\\ \\ \displaystyle +\frac{vv_{2}}{1+v}{\tilde D}^{a}v{\tilde D}_{a}u-\frac{2v(v_{2})^{2}}{1+v}.\end{array}
\end{equation}
D'autre part, la d\'efinition des composantes $A^{ab}$ montre que
\begin{equation}
G^{ab}=\frac{A^{ab}+{\tilde D}^{a}u{\tilde D}^{b}u}{1+v}.
\end{equation}
De ce fait, la somme des deux derniers termes du membre de droite de $(4.8)$ s'\'ecrit sous la forme:
$$\begin{array}{c}\displaystyle {\tilde D}_{c}\ ^{a}u{\tilde D}^{c}u{\tilde D}^{b}u{\tilde D}_{ab}u+{\tilde D}^{a}u{\tilde D}^{c}u{\tilde D}_{c}\ ^{b}u{\tilde D}_{ab}u=\frac{2\Big( {\tilde D}^{a}u{\tilde D}^{b}u{\tilde D}_{ab}u\Big)^{2}}{1+v}\\ \\\displaystyle +\frac{A^{ad}{\tilde D}_{cd}u{\tilde D}^{c}u{\tilde D}^{b}u{\tilde D}_{ab}u}{1+v}+\frac{A^{ad}{\tilde D}_{cd}u{\tilde D}^{c}u{\tilde D}^{b}u{\tilde D}_{ba}u}{1+v}.\end{array}$$
Or, un simple calcul nous permet d'\'etablir la commutation suivante :
\begin{equation}
D_{cd}u=D_{dc}u-T^{e}_{cd}D_{e}u.
\end{equation}
Il en d\'ecoule que
$$\begin{array}{c}\displaystyle {\tilde D}_{c}\ ^{a}u{\tilde D}^{c}u{\tilde D}^{b}u{\tilde D}_{ab}u+{\tilde D}^{a}u{\tilde D}^{c}u{\tilde D}_{c}\ ^{b}u{\tilde D}_{ab}u=\frac{2\left( {\tilde D}^{a}u{\tilde D}^{b}u{\tilde D}_{ab}u\right) ^{2}}{1+v}\\ \\ \displaystyle +\frac{2A^{ab}({\tilde D}_{ac}u{\tilde D}^{c}u)({\tilde D}_{bd}u}{\tilde D}^{d}u){1+v}+E_{3}+E_{4},\end{array}$$\rm
o\`u l'on a not\'e 
$$E_{3}=-\frac{3}{1+v}A^{ab}e^{(\mu_{a}+\mu _{c})u}T^{e}_{ca}D_{e}u{\tilde D}^{c}u{\tilde D}^{d}u{\tilde D}_{bd}u$$
et
\begin{equation}
E_{4}=\frac{1}{1+v}A^{ab}e^{(\mu _{a}+\mu _{b}+\mu _{c}+\mu _{d})u}T^{e}_{cb}T^{f}_{da}D_{e}uD_{f}u{\tilde D}^{c}u{\tilde D}^{d}u.
\end{equation}
Ainsi, tenons compte de $(4.9)$, on v\'erifie que
\begin{equation}
\begin{array}{c}\displaystyle {\tilde D}_{c}\ ^{a}u{\tilde D}^{c}u{\tilde D}^{b}u{\tilde D}_{ab}u+{\tilde D}^{a}u{\tilde D}^{c}u{\tilde D}_{c}\ ^{b}u{\tilde D}_{ab}u=\frac{2\Big( {\tilde D}^{a}u{\tilde D}^{b}u{\tilde D}_{ab}u\Big) ^{2}}{1+v}\\ \\ \displaystyle +\frac{A^{ab}{\tilde D}_{a}v{\tilde D}_{b}v}{2(1+v)}+\frac{2(v_{2})^{2}v}{1+v}-\frac{2v_{2}}{1+v}A^{ab}{\tilde D}_{a}v{\tilde D}_{b}u+E_{3}+E_{4}.\end{array}
\end{equation}
D'autre part, utilisons \`a nouveau $(4.9)$, on obtient la forme suivante pour le terme $E_{3}$:
$$\begin{array}{c}\displaystyle E_{3}=-\frac{3}{2(1+v)}A^{ab}e^{(\mu _{a}+\mu _{c})u}T^{e}_{ca}D_{e}u{\tilde D}^{c}u{\tilde D}_{b}v\\ \\ \displaystyle +\frac{3v_{2}}{1+v}A^{ab}e^{(\mu _{a}+\mu _{c})u}T^{e}_{ca}D_{e}u{\tilde D}^{c}u{\tilde D}_{b}u.\end{array}$$
Une forme, qui compte tenu de $(4.5)$ et de la d\'efinition des composantes $A^{ab}$, permet de v\'erifier que
\begin{equation}
E_{3}=\Big(-\frac{3l}{2}+\frac{3v_{2}}{1+v}\Big)e^{(\mu _{a}+\mu _{c})u}T^{e}_{ca}D_{e}u{\tilde D}^{c}u{\tilde D}^{a}u.
\end{equation}

\vskip0mm

Reportons $(4.10)$ et $(4.14)$ dans $(4.8)$, tenons compte \`a nouveau de la d\'efinition des composantes $A^{ab}$, on obtient :
$$\begin{array}{c}\displaystyle A^{ab}{\tilde D}_{cab}u{\tilde D}^{c}u=\frac{A^{ab}{\tilde D}_{a}v{\tilde D}_{b}v}{2(1+v)}+{\tilde D}_{c}F{\tilde D}^{c}u-\frac{F{\tilde
D}_{a}v{\tilde D}^{a}u}{1+v}-\frac{v_{2}(2+v)}{1+v}{\tilde D}_{a}v{\tilde D}^{a}u\\ \\ \displaystyle -v(1+v)(\mu _{a}+\mu _{b})G^{ab}{\tilde D}_{ab}u+4v\mu _{a}{\tilde D}^{a}u{\tilde D}^{b}u{\tilde D}_{ab}u+\frac{4v(v_{2})^{2}}{1+v}+E_{3}+E_{4}.\end{array}$$
Une \'egalit\'e qui, compte tenu de $(4.9)$, s'\'ecrit sous la forme :
$$\begin{array}{c}\displaystyle A^{ab}{\tilde D}_{cab}u{\tilde D}^{c}u=\frac{A^{ab}{\tilde D}_{a}v{\tilde D}_{b}v}{2(1+v)}+{\tilde D}_{c}F{\tilde D}^{c}u-\frac{F{\tilde
D}_{a}v{\tilde D}^{a}u}{1+v}-\frac{v_{2}(2+v)}{1+v}{\tilde
D}_{a}v{\tilde D}^{a}u\\ \\ \displaystyle -v(1+v)(\mu _{a}+\mu
_{b})G^{ab}{\tilde D}_{ab}u+2v\mu _{a}{\tilde D}^{a}v{\tilde
D}_{a}u-\frac{4(vv_{2})^{2}}{1+v}+E_{3}+E_{4}.\end{array}$$
Multiplions cette relation par 2 et reportons dans $(4.7)$, l'\'egalit\'e qui
en r\'esulte s'\'ecrit sous la forme :
\begin{equation}
\begin{array}{c}\displaystyle A^{ab}{\tilde D}_{ab}v-\frac{A^{ab}{\tilde D}_{a}v{\tilde D}_{b}v}{(1+v)}=2{\tilde D}_{a}F{\tilde
D}^{a}u-2\frac{F{\tilde D}_{a}v{\tilde D}^{a}u}{1+v}-\frac{2v_{2}(2+v)}{1+v}{\tilde D}_{a}v{\tilde D}^{a}u\\ \\ \displaystyle -2v(1+v)(\mu _{a}+\mu _{b})G^{ab}{\tilde
D}_{ab}u+8A^{ab}\mu _{c}{\tilde D}_{ac}u{\tilde D}^{c}u{\tilde D}_{b}u+4v\mu _{a}{\tilde D}_{a}v{\tilde D}^{a}u\\ \\ \displaystyle +2v_{2}F+4vv_{2}-\frac{8(vv_{2})^{2}}{1+v}+2A^{ab}{\tilde D}^{c}_{a}u{\tilde D}_{bc}u+4E_{1}+E_{2}+2E_{3}+2E_{4}.\end{array}
\end{equation} 
La d\'efinition des composantes $A^{ab}$ montre, \`a nouveau, que 
$$2\mu _{c}A^{ab}{\tilde D}_{ac}u{\tilde D}^{c}u{\tilde D}_{b}u=(\mu _{a}+\mu _{b}){\tilde D}^{a}u{\tilde D}^{b}u{\tilde D}_{ab}u.$$
D'autre part, d'apr\`es $(4.9)$, on obtient :
\begin{equation}
(\mu _{a}+\mu _{b}){\tilde D}^{a}u{\tilde D}^{b}u{\tilde D}_{ab}u=2\mu _{a}{\tilde
D}^{a}u{\tilde D}^{b}u{\tilde D}_{ab}u=\mu _{a}{\tilde D}^{a}u{\tilde D}_{a}v-2(v_{2})^{2}.
\end{equation}
Reportons dans $(4.16)$, on obtient :
$$\begin{array}{c}\displaystyle A^{ab}{\tilde D}_{ab}v-\frac{A^{ab}{\tilde D}_{a}v{\tilde D}_{b}v}{(1+v)}=2{\tilde D}_{c}F{\tilde D}^{c}u-2\frac{F{\tilde D}_{a}v{\tilde D}^{a}u}{1+v}-\frac{2v_{2}(2+v)}{1+v}{\tilde D}_{a}v{\tilde D}^{a}u\\ \\ \displaystyle -2v(1+v)(\mu _{a}+\mu _{b})G^{ab}{\tilde D}_{ab}u+4v\mu _{a}{\tilde D}_{a}v{\tilde D}^{a}u+4\mu _{a}{\tilde D}^{a}u{\tilde D}_{a}v+2v_{2}F\\ \\ \displaystyle +4v_{2}(v_{1}-v_{2})-\frac{8(vv_{2})^{2}}{1+v}+2A^{ab}{\tilde D}^{c}_{a}u{\tilde D}_{bc}u+4E_{1}+E_{2}+2E_{3}+2E_{4}.\end{array}$$
Supposons que
\begin{equation}
v(\xi )\geq 1.
\end{equation}
Tenons compte de $(4.2)$ et $(4.5)$, on voit que
\begin{equation}
\begin{array}{c}\displaystyle A^{ab}{\tilde D}_{ab}v-\frac{A^{ab}{\tilde D}_{a}v{\tilde
D}_{b}v}{(1+v)}\geq 2{\tilde D}_{c}F{\tilde D}^{c}u-2lvF+2l(1+v)vv_{2}-8v(v_{2})^{2}\\ \\ \displaystyle +2(m+l)vv_{2}+4lv_{2}+2\alpha (1+v)v^{\frac{3}{2}}-2(1+v)^{\frac{3}{2}}v_{2}e^{u}K(e^{u}\xi
) \\ \\ \displaystyle +2v_{1}v_{2}-\frac{8(v_{2})^{2}}{1+v}
+2A^{ab}{\tilde D}^{c}_{a}u{\tilde D}_{bc}u+4E_{1}+E_{2}+2E_{3}+2E_{4}.\end{array}
\end{equation}
Remarquons que l'\'egalit\'e $(4.17)$ implique aussi la suivante : 
$${\tilde D}_{c}(v_{2}){\tilde D}^{c}u=\mu _{a}{\tilde
D}_{a}v{\tilde D}^{a}u+2v_{2}(v-v_{2})$$
de sorte que le d\'eveloppement de ${\tilde D}_{c}F{\tilde D}^{c}u$ donne
$$\begin{array}{c}\displaystyle 2{\tilde D}_{c}F{\tilde D}^{c}u=-4v_{2}(v-v_{1})-2\mu
_{a}{\tilde D}_{a}v{\tilde D}^{a}u+2(m-1){\tilde D}_{a}v{\tilde D}^{a}u\\ \\
\displaystyle -\sqrt{1+v}e^{u}\Big[3K(e^{u}\xi ){\tilde D}_{a}v{\tilde
D}^{a}u+2(1+v)({\tilde D}_{c}K) (e^{u}\xi ){\tilde D}^{c}u\Big]\\ \\
\displaystyle -2v(1+v)^{\frac{3}{2}}e^{u}\frac{\partial [\rho K(\rho \xi )]}{\partial \rho }(e^{u}\xi ).\end{array}$$
Utilisons $(4.5)$, on montre que
$$\begin{array}{c}\displaystyle 2{\tilde D}_{c}F{\tilde D}^{c}u-2lvF=l(1+v)F-l(m-1)(1+v)^{2}+l(1+v)v_{2}\\ \\ \displaystyle -2(1+v)^{\frac{3}{2}}e^{u}\left[v\frac{\partial [\rho K(\rho \xi )]}{\partial \rho }(e^{u}\xi )+({\tilde D}_{c}K)(e^{u}\xi ){\tilde D}^{c}u\right]\\ \\ \displaystyle -4(v_{2})^{2}-2lv_{2}+l(1+v)^{\frac{3}{2}}e^{u}K(e^{u}\xi ).\end{array}$$
Reportons cette \'egalit\'e dans $(4.19)$, la relation qui en r\'esulte se transforme comme suit : 
$$\begin{array}{c}\displaystyle A^{ab}{\tilde D}_{ab}v-\frac{A^{ab}{\tilde
D}_{a}v{\tilde D}_{b}v}{(1+v)}-l(1+v)F\geq
l(1+v)^{\frac{3}{2}}e^{u}K(e^{u}\xi )+2(v_{2})^{2}\\ \\ \displaystyle 
+\frac{8v_{2}(1+v_{1})}{1+v}+l(1+v)v_{2}+2(m+5)vv_{2}+8(1+v)v_{1}v_{2}\cr \cr
\displaystyle -l(m-1)(1+v)^{2}-2(1+v)^{\frac{3}{2}}e^{u}({\tilde D}_{a}K)(e^{u}\xi ){\tilde D}^{a}u\\ \\ \displaystyle +2v(1+v)\Big[\alpha \sqrt{v}-\sqrt{1+v}e^{u}\frac{\partial [\rho K(\rho \xi )]}{\partial \rho }(e^{u}\xi )\Big]\\ \\ \displaystyle +2(1+v)^{\frac{3}{2}}v_{2}\left[(l-4)\sqrt{1+v}-e^{u}K(e^{u}\xi)\right]\\ \\ \displaystyle +2A^{ab}{\tilde D}^{c}_{a}u{\tilde D}_{bc}u+4E_{1}+E_{2}+2E_{3}+2E_{4}.\end{array}$$
Supposons que $\displaystyle l\geq l_{0}=4+r_{2}\max _{\Sigma _{r_{1},r_{2}}}K$, tenons compte du choix $(4.3)$ de $\alpha $ et minorons par $0$ les termes positifs, l'\'egalit\'e pr\'ec\'edente permet d'\'ecrire
\begin{equation}
\begin{array}{c}\displaystyle A^{ab}{\tilde D}_{ab}v-\frac{A^{ab}{\tilde D}_{a}v{\tilde
D}_{b}v}{(1+v)}-l(1+v)F\geq -l(m-1)(1+v)^{2}+2A^{ab}{\tilde
D}^{c}_{a}u{\tilde D}_{bc}u\\ \\ \displaystyle -2(1+v)^{2}e^{u}\Vert {\tilde D}K\Vert _{\infty }+4E_{1}+E_{2}+2E_{3}+2E_{4}.\end{array}
\end{equation}
Or, $\vert E_{2}\vert \leq C_{1}v^{2}$, o\`u $C_{1}$ est fonction de $\Vert {\cal
T}\Vert _{\infty }$, $\Vert \nabla {\cal T}\Vert _{\infty }$, $\Vert {\tilde
{\cal R}}\Vert _{\infty }$ et $\Vert u\Vert _{\infty }$ et puisque, d'apr\`es $(4.13)$, on sait que $\displaystyle \vert E_{4}\vert \leq C'_{1}v^{2}$ et compte tenu de $(4.15)$, on a : $\displaystyle \vert E_{3}\vert \leq C_{2}lv^{\frac{3}{2}}$, o\`u les constantes $C'_{1}$ et $C_{2}$ ne sont fonction que de $\Vert {\cal T}\Vert _{\infty }$ et $\Vert u\Vert _{\infty }$, l'in\'egalit\'e $(4.20)$ implique l'existence d'une constante positive $C_{3}$ telle que
$$\begin{array}{c}\displaystyle A^{ab}{\tilde D}_{ab}v-\frac{1}{(1+v)}A^{ab}{\tilde D}_{a}v{\tilde D}_{b}v+l(1+v)F\geq -2C_{3}l(1+v)^{2}\\ \\ \displaystyle +2A^{ab}{\tilde D}^{c}_{a}u{\tilde D}_{bc}u+4E_{1}.\end{array}$$
Reportons cette in\'egalit\'e dans $(4.6)$, on oboutit \`a la suivante :
\begin{equation}
-C_{3}l(1+v)^{2}+A^{ab}{\tilde D}^{c}_{a}u{\tilde D}_{bc}u+2E_{1}\leq 0.
\end{equation}

\vskip2mm

Maintenant, on d\'eveloppe le carr\'e suivant :
$$\begin{array}{c}\displaystyle K=G^{ab}G^{cd}\Big[{\tilde D}_{ad}u-\epsilon {\tilde
D}_{a}u{\tilde D}^{e}u{\tilde D}_{ed}u-e^{(\mu _{a}+\mu _{e})u}T^{l}_{ae}{\tilde D}^{e}uG_{ld}\Big]\\ \\ \displaystyle \Big[{\tilde D}_{bc}u-\epsilon {\tilde D}_{b}u{\tilde D}^{e}u{\tilde D}_{ec}u-e^{(\mu _{b}+\mu _{e})u}T^{l}_{be}{\tilde D}^{e}uG_{lc}\Big]\end{array}$$
avec $\displaystyle \epsilon =\frac{1}{v}$ pour voir que
$$\begin{array}{c}\displaystyle (1+v)K=A^{ab}{\tilde D}^{c}_{a}u{\tilde
D}_{bc}u+(1+v)e^{(\mu _{a}+\mu _{e}+\mu _{b}+\mu _{f})u}G^{ab}G_{cd}T^{c}_{ae}{\tilde
D}^{e}uT^{d}_{bf}{\tilde D}^{f}u\\ \\ \displaystyle +2E_{1}-v^{-1}{\tilde D}^{a}u{\tilde D}^{b}u{\tilde D}_{ac}u{\tilde D}_{b}^{c}u+2v^{-1}e^{(\mu _{b}+\mu
_{e})u}T^{d}_{be}{\tilde D}^{e}u{\tilde D}^{a}u{\tilde D}^{b}u{\tilde
D}_{ad}u.\end{array}$$
Il existe, alors, une constante positive $C_{4}$ fonction de $\Vert u\Vert
_{\infty }$ et $\Vert {\cal T}\Vert _{\infty }$, telle que
\begin{equation}
\begin{array}{c}\displaystyle A^{ab}{\tilde D}^{c}_{a}u{\tilde
D}_{bc}u+2E_{1}\geq -C_{4}(1+v)^{2}+v^{-1}{\tilde D}^{a}u{\tilde D}^{b}u{\tilde
D}_{ac}u{\tilde D}_{b}^{\ c}u\\ \\ \displaystyle -2v^{-1}e^{(\mu _{b}+\mu
_{c})u}{\tilde D}^{a}u{\tilde D}^{b}u{\tilde D}^{d}_{\ a}uT^{e}_{bc}{\tilde
D}_{e}u.\end{array}
\end{equation}
Utilisons $(4.12)$, on montre que
$$\begin{array}{c}\displaystyle {\tilde D}^{a}u{\tilde D}^{b}u{\tilde D}_{ac}u{\tilde
D}_{b}^{\ c}u={\tilde D}^{a}u{\tilde D}^{b}u{\tilde D}_{ca}u{\tilde
D}^{c}_{\ b}u-2{\tilde D}^{a}u{\tilde D}^{b}u{\tilde D}^{c}_{\ a}ue^{(\mu
_{b}+\mu _{c})u}T^{e}_{bc}{\tilde D}_{e}u\\ \\ \displaystyle +G^{cd}{\tilde
D}^{a}u{\tilde D}^{b}ue^{(\mu _{a}+\mu _{b}+\mu _{c}+\mu _{d})u}T^{e}_{ac}{\tilde
D}_{e}uT^{f}_{bd}{\tilde D}_{f}u\end{array}$$
et,
donc, d'apr\`es $(4.9)$, on aura :
$$\begin{array}{c}\displaystyle {\tilde D}^{a}u{\tilde D}^{b}u{\tilde D}_{ac}u{\tilde
D}_{b}^{\ c}u=\frac{1}{4}{\tilde D}^{a}v{\tilde D}_{a}v-v_{2}{\tilde
D}^{a}v{\tilde D}_{a}u+v(v_{2})^{2} \\ \\ \displaystyle -({\tilde
D}^{b}u{\tilde D}^{c}v-2v_{2}{\tilde D}^{b}u{\tilde D}^{c}u)e^{(\mu
_{b}+\mu _{c})u}T^{e}_{bc}{\tilde D}_{e}u\\ \\ \displaystyle
+G^{cd}{\tilde D}^{a}u{\tilde D}^{b}ue^{(\mu _{a}+\mu _{b}+\mu _{c}+\mu
_{d})u}T^{e}_{ac}{\tilde D}_{e}uT^{f}_{bd}{\tilde D}_{f}u.\end{array}$$
Ainsi, tenons compte de $(4.5)$, il existe une constante positive $C_{5}$ fonction de $\Vert u\Vert _{\infty }$ et $\Vert {\cal T}\Vert _{\infty }$, telle que
\begin{equation}
{\tilde D}^{a}u{\tilde D}^{b}u{\tilde D}_{ac}u{\tilde D}_{b}^{\ c}u\geq
\frac{l^{2}}{4}(1+v)^{2}v-l(1+v)vv_{2}-C_{5}l(1+v)v_{2}\sqrt{v}.
\end{equation}
De m\^eme, on v\'erifie que
$$\begin{array}{c}\displaystyle 2{\tilde D}^{a}u{\tilde D}^{b}u{\tilde D}_{ad}uT^{d}_{be}{\tilde D}^{e}u=T^{d}_{be}{\tilde D}^{e}u{\tilde D}^{b}u{\tilde
D}_{d}v-2v_{2}T^{d}_{be}{\tilde D}^{b}u{\tilde D}^{e}u{\tilde D}_{d}u\\ \\
\displaystyle -2e^{(\mu _{a}+\mu _{d})u}T^{d}_{be}T^{f}_{ad}{\tilde
D}^{a}u{\tilde D}^{b}u{\tilde D}^{e}u{\tilde D}_{f}u\end{array}$$
et donc $(4.5)$ implique l'existence d'une constante positive $C_{6}$ fonction de
$\Vert u\Vert _{\infty }$ et $\Vert {\cal T}\Vert _{\infty }$, telle que
\begin{equation}
\vert 2e^{(\mu _{b}+\mu _{c})u}{\tilde D}^{a}u{\tilde D}^{b}u{\tilde D}^{d}_{\
a}uT^{e}_{bc}{\tilde D}_{e}u\vert \leq C_{6}l(1+v)v_{2}\sqrt{v}.
\end{equation}
Reportons $(4.23)$ et $(4.24)$ dans $(4.22)$, compte tenu du fait que $v_{2}\leq 1+v$, on en d\'eduit que
$$A^{ab}{\tilde D}^{c}_{a}u{\tilde D}_{bc}u+2E_{1}\geq
(4^{-1}l^{2}-l-C_{4})(1+v)^{2}-lC_{7}(1+v)^{\frac{3}{2}},$$
o\`u $C_{7}=C_{5}+C_{6}$. Une in\'egalit\'e qui, compte tenu de $(4.21)$, implique 
$$(4^{-1}l^{2}-lC_{3}-l-C_{4})\sqrt{1+v}-lC_{7}\leq 0.$$
De sorte que, pour $l$ assez grand, on aura : $1+v(\xi )\leq C_{8}=(lC_{7})^{2}$.
Ainsi, compte tenu de $(4.18)$, $1+v\leq 1+C_{8}$. La d\'efinition $(4.4)$ de la fonctionnelle $\Gamma $ montre que partout dans $\Sigma $, on a: 
$$v\leq (1+C_8)\left(\frac{r_2}{r_1}\right) ^l.$$
Le lemme est prouv\'e.

\vskip6mm

\begin{lem}Soient \it $K\in \mathscr{C}^{\infty }(E_{*})$ une fonction partout strictement positive v\'erifiant l'hypoth\`ese $(1.1)$ du th\'eor\`eme $2$. Alors toute
solution $u\in C^{3}(\Sigma )$ d'une \'equation de la forme
\begin{equation}
\sum _{n+1\leq \alpha ,\beta \leq n+m-1}B^{\alpha \beta }(u)D_{\alpha \beta
}u=(m-1)(1+v_{1})-(1+v_{1})^{\frac{3}{2}}e^{u}K(e^{u}\xi ),
\end{equation}
o\`u $B^{\alpha \beta }=(1+v_{1})G^{\alpha \beta }-D^{\alpha }uD^{\beta }u$,
v\'erifie $r_{1}\leq e^{u}\leq r_{2}$. Si en outre $K$ v\'erifie l'hypoth\`ese de
monotonicit\'e $(1.2)$ du th\'eor\`eme $2$, il existe alors une constante positive
$C_{0}$ ne d\'ependant que de $n$, $r_{1}$, $r_{2}$, $\displaystyle \max
_{\Sigma _{r_{1},r_{2}}}K$ et $\displaystyle \Vert K\Vert _{\mathscr{C}^{1}(\Sigma
_{r_{1},r_{2}})}$ telle que $\displaystyle v_{1}\leq C_{0}$.
\end{lem}

\vskip5mm

\noindent\textit{D\'emonstration. }Soit $u\in \mathscr{C}^{3}(\Sigma )$ une solution de $(4.25)$. Soit $\xi \in \Sigma$ un point o\`u $u$ atteint son maximum. Si $u(\xi )>log(r_{2})$, l'hypoth\`ese de croissance $(1.1)$ faite sur $K$ implique qu'au point $\xi $, on aura 
$$0\geq {\tilde D}_{\alpha }^{\alpha }u=(m-1)-e^{u}K(e^{u}\xi )>0$$
ce qui constitue une contradiction. La minoration $u\geq log(r_{1})$ s'obtient par analogie en consid\'erant un point o\`u $u$ atteint son minimum.

\vskip3mm

A pr\'esent, on d\'erive une fois l'\'equation $(4.25)$ suivant une direction verticale, il vient
$$B^{\alpha \beta }D_{\lambda \alpha \beta }u+D_{\lambda }B^{\alpha \beta }D_{\alpha
\beta }u=D_{\lambda }F,$$
o\`u $F$ d\'esigne le second membre de $(4.25)$. D\'eveloppons $D_{\lambda }B^{\alpha \beta }$ et saturons par $D^{\lambda }u$ , on obtient,
\begin{equation}
B^{\alpha \beta }D_{\lambda \alpha \beta }uD^{\lambda }u+\bigtriangleup
_{1}uD_{\lambda }v_{1}D^{\lambda }u-\frac{1}{2}D_{\lambda }v_{1}D^{\lambda }{\tilde
v}=D_{\lambda }FD^{\lambda }u,
\end{equation}
o\`u l'on a not\'e $\bigtriangleup _{1}u=G^{\alpha \beta }D_{\alpha \beta }u$. Ce terme
peut \^etre exprim\'e \`a partir de $(4.25)$ de la mani\`ere suivante:
$$G^{\alpha \beta }D_{\alpha \beta }u=\frac{F+D^{\alpha }uD^{\beta }uD_{\alpha \beta
}u}{1+v_{1}}$$
et puisque $\displaystyle D^{\beta }uD_{\alpha \beta }u=\frac{1}{2}D_{\alpha }v_{1}$, on obtient
$$G^{\alpha \beta }D_{\alpha \beta }u=\frac{F}{1+v_1}+\frac{1}{2}\frac{D^{\alpha
}uD_{\alpha }v_1}{1+v_1}.$$
Par suite
$$\bigtriangleup _{1}uD_{\lambda }v_{1}D^{\lambda }u=\frac{F}{1+v_1}D_{\lambda
}v_{1}D^{\lambda }u+\frac{1}{2}\frac{D^{\alpha }uD^{\lambda }uD_{\alpha
}v_{1}D_{\lambda }v_1}{1+v_1}.$$
Mais, d'apr\`es l'expression des composantes $B^{\alpha \beta }$ qui implique que
\begin{equation}
D^{\alpha }uD^{\beta }u=(1+v_{1})G^{\alpha \beta}-B^{\alpha \beta },
\end{equation}
on peut \'ecrire
\begin{equation}
\bigtriangleup _{1}uD_{\lambda }v_1D^{\lambda}u=\frac{F}{1+v_1}D_{\lambda
}v_1D^{\lambda }u+\frac{1}{2}D_{\lambda }v_{1}D^{\lambda }v_1-\frac{1}{2}\frac{B^{\alpha \beta}D_{\alpha}v_1D_{\beta }v_1}{1+ v_1}.
\end{equation}
D'autre part, et puisque $\displaystyle D_{\lambda \alpha \beta }u=D_{\alpha \beta
\lambda }u+{\tilde R}^{A}_{\beta \alpha \lambda }D_{A}u$, les expressions $(2.15)$ et $(2.16)$ donnant les composantes de la forme ${\tilde R}^{A}_{\beta \alpha \lambda }$ du tenseur de courbure, impliquent que
\begin{equation}
B^{\alpha \beta }D_{\lambda \alpha \beta }uD^{\lambda }u=B^{\alpha \beta }D_{\alpha
\beta \lambda }uD^{\lambda }u-(m-2)v_{1}(1+ v_{1}).
\end{equation}
Reportons $(4.28)$ et $(4.29)$ dans $(4.26)$ et multiplions par $2$, la relation qui en r\'esulte s'\'ecrit sous la forme :
\begin{equation}
\begin{array}{c}\displaystyle 2B^{\alpha \beta }D_{\alpha \beta
\lambda}uD^{\lambda }u-\frac{B^{\alpha \beta }D_{\alpha }v_1D_{\beta
}v_1}{1+v_1}=2(m-2)v_{1}(1+v_{1})\\ \displaystyle +2D_{\lambda
}FD^{\lambda }u-\frac{2FD_{\lambda }v_1D^{\lambda }u}{1+v_1}.\end{array}
\end{equation}

\vskip0mm

Maintenant, on consid\`ere la fonctionnelle
\begin{equation}
\Gamma (u)=(1+v_{1})e^{lu},
\end{equation}
o\`u $l$ est un r\'eel strictement positif fix\'e ult\'erieurement. En un point $\xi \in \Sigma $ o\`u $\Gamma (u) $ atteint son maximum, on a :
\begin{equation}
_frac{D_{\lambda }\Gamma }{\Gamma }=\frac{D_{\lambda }v_{1}}{1+v_{1}}+lD_{\lambda
}u=0
\end{equation}
et $B^{\alpha \beta }D_{\alpha \beta }\Gamma \leq 0$, c'est-\`a-dire, tenons
compte du fait que $u$ v\'erifie l'\'equation $(4.25)$,
$$\displaystyle B^{\alpha \beta }D_{\alpha \beta }v_1-\frac{B^{\alpha
\beta }D_{\alpha }v_1D_{\beta }v_1}{1+v_1}+l(1+v_1)F\leq 0.$$
Or,
$$B^{\alpha \beta }D_{\alpha \beta }v_1=2B^{\alpha \beta }D_{\alpha \beta \lambda }uD^{\lambda }u+2B^{\alpha \beta }D_{\alpha \lambda }uD_{\beta }^{\lambda }u.$$
Donc, d'apr\`es $(4.30)$ et puisque $(m-2)v_1(1+v_1)\geq 0$, au maximum de $\Gamma$, on doit avoir 
$$2B^{\alpha \beta }D_{\alpha \lambda }uD_{\beta }^{\lambda
}u+l(1+v_1)F+2D_{\lambda }FD^{\lambda }u-\frac{2FD_{\lambda }v_1D^{\lambda }u}{1+
v_1}\leq
0$$
c'est-\`a-dire, tenons compte de $(4.32)$,
\begin{equation}
2B^{\alpha \beta }D_{\alpha \lambda }uD_{\beta
}^{\lambda}u+lF+3lv_{1}F+2D_{\lambda }FD^{\lambda }u\leq 0.
\end{equation}
D\'eveloppons $D_{\lambda }FD^{\lambda }u$. Compte tenu de $(4.32)$, on v\'erifie que
$$\begin{array}{c}\displaystyle 2D_{\lambda
}FD^{\lambda
}u+3lv_1F=-2(1+v_1)^{\frac{3}{2}}e^{u}(D_{\lambda
}K)(e^{u}\xi )D^{\lambda }u\\ \\ \displaystyle
+l(m-1)(1+v_1)v_1-2v_1(1+v_1)^{\frac{3}{2}}e^{u}\frac{\partial [\rho
K(\rho\xi )]}{\partial \rho}(e^{u}\xi ).\end{array}$$
L'hypoth\`ese de croissance $(1.2)$ faite sur $K$ implique donc l'existence d'une
constante positive $C_{1}$, fonction de $\displaystyle \Vert K\Vert _{\mathscr{C}^{1}(\Sigma _{r_1,r_2})}$, telle que
$$2D_{\lambda }FD^{\lambda }u+3lv_1F\geq (l-C_{1})(m-1)(1+v_1)v_1.$$
Reportons dans $(4.33)$, on obtient
$$2B^{\alpha \beta }D_{\alpha \lambda }uD_{\beta }^{\lambda
}u+lF+(l-C_{1})(m-1)(1+v_1)v_1\leq 0$$
et par suite, minorons par $0$ les termes positifs, on aboutit \`a l'in\'egalit\'e
suivante:
\begin{equation}
(l-C_{1})(m-1)(1+v_1)v_1-l(1+v_1)^{\frac{3}{2}}e^{u}K(e^{u}\xi
)\leq 0.
\end{equation}
Posons $l=1+C_{1}$. Si
\begin{equation}
v_{1}(\xi )\geq 1,
\end{equation}
l'in\'egalit\'e $(4.34)$ combin\'ee avec le fait que $m\geq 2$ montre que
$$\sqrt{v_1}\leq 3\sqrt{2}(1+C_1)e^{u}K(e^{u}\xi
)\leq 3\sqrt{2}r_2(1+C_1)K(e^{u}\xi ).$$
En cons\'equence, compte tenu de $(4.35)$, on voit que
$$v_1(\xi )\leq C_2=1+18\left[r_2(1+C_1)\max _{\Sigma
_{r_1,r_2}}K\right]^{2}.$$
Enfin, usons de la d\'efinition $(4.31)$ de la fonctionnelle $\Gamma $, on v\'erifie que, partout dans $\Sigma $, on a :
$$v_1\leq 1+v_1\leq C_2\left(\frac{r_2}{r_1}\right)^{l}.$$
Le lemme est prouv\'e.

\vskip6mm

\section{Preuve des r\'esultats}
\subsection{Preuve du th\'eor\`eme 1}

\vskip4mm

Les calculs de la premi\`ere section montrent que la courbure moyenne de $\Sigma _{r}$ est
$${\cal {M}}_{_{\Sigma _r}}=-G^{ab}G(D_{e_{a}}e_{b},\nu )=\frac{m-1}{r}.$$
D\'efinissons ${\cal Y}$ sur $\Sigma $ en posant 
$${\cal Y}(\xi )=r\xi \mbox{ pour tout }\xi \in \Sigma .$$
Le graphe ${\cal Y}$ n'est autre que $\Sigma _{r}$ et s'il existe un r\'eel $r>0$ tel que 
$$K(r\xi )=(m-1)r^{-1},\ \mbox{quel que soit }\xi \in \Sigma,$$
alors 
$${\cal {M}}_{_{\cal Y}}(r\xi)=K(r\xi ),\mbox{ pour tout }\xi \in \Sigma .$$
Ceci montre la suffisance des hypoth\`eses faites sur $K$ dans le th\'eor\`eme 1.

\vskip3mm

Supposons qu'il existe un graphe radial $\cal {Y}$ sur $\Sigma $ de classe $\mathscr{C}^{\infty }$ dont la courbure moyenne est \'egale \`a $K$. Ceci est \'equivaut \`a dire
qu'il existe une fonction $u\in C^{\infty }(\Sigma )$ 
v\'erifiant l'\'equation $(3.4)$, mise en \'evidence \`a la troisi\`eme section. Notons
$$r_1=\min _{\Sigma }e^{u}\mbox{ et }r_2=\max _{\Sigma }e^{u}.$$
On se place successivement en un point o\`u $u$ atteint son minimum et en un point o\`u $u$ atteint son maximum, on voit que  
$$(m-1)(r_1)^{-1}\geq K(r_1)\mbox{ et }(m-1)(r_2)^{-1}\leq K(r_2).$$
Le r\'esultat devient une cons\'equence de la continuit\'e de $K$ via le th\'eor\`eme
des valeurs interm\'ediaires.

\vskip6mm

\subsection{Preuve du th\'eor\`eme 2}

\vskip4mm

Pour tout entier naturel $k$ et tout r\'eel $\alpha $ compris entre $0$ et $1$, on notera par $A^{k,\alpha }(\Sigma )$ l'ensemble des fonctions $u\in \mathscr{C}^{k,\alpha }(\Sigma )$ dont la composante horizontale du gradient est identiquement
nulle muni de la norme $\mathscr{C}^{k,\alpha }(\Sigma )$. 

\vskip4mm

La preuve du th\'eor\`eme repose sur le lemme suivant.

\vskip5mm

\begin{lem}Soit $w\in A^{k+1,\alpha }(\Sigma )$. Pour tout $u\in \mathscr{C}^{k+2,\alpha }(\Sigma )$, on pose
$$L[w](u)=\sum _{n+1\leq \alpha ,\beta \leq n+m-1}B^{\alpha \beta }(w)D_{\alpha \beta
}u-u,$$
o\`u $B^{\alpha \beta }(w)=(1+v_{1}(w))G^{\alpha \beta }-D^{\alpha }wD^{\beta }w$ et
$v_{1}(w)=D_{\alpha }wD^{\alpha }w$ est le carr\'e de la norme de la composante verticale du gradient de $w$. L'op\'erateur $L[w]$ ainsi d\'efini r\'ealise un isomorphisme de $A^{k+2,\alpha }(\Sigma )$ sur $A^{k,\alpha }(\Sigma )$.
\end{lem}

\vskip5mm

\noindent\textit{D\'emonstration. }Tout d'abord, on applique la m\'ethode de continuit\'e, combin\'ee avecles in\'egalit\'es de Schauder, pour montrer que pour toute fonction $w\in C^{k+1,\alpha }(\Sigma )$, l'op\'erateur
elliptique ${\tilde L}[w]$ d\'efini, pour $u\in C^{k+2,\alpha }(\Sigma )$, par
$${\tilde L}[w]u=G^{ij}D_{ij}u+B^{\alpha \beta }(w)D_{\alpha \beta }u-u$$ 
est un isomorphisme de $C^{k+2,\alpha }(\Sigma )$ sur $C^{k,\alpha
}(\Sigma )$. On en d\'eduit que pour toute fonction $w\in A^{k+1,\alpha }(\Sigma )$,
l'op\'erateur $L[w]$ est un isomorphisme de $A^{k+2,\alpha }(\Sigma )$ sur $A^{k,\alpha
}(\Sigma )$. En effet, le principe du maximum implique que le noyau de celui-ci est
trivial et pour toute fonction $z\in A^{k,\alpha }(\Sigma )$, il existe $u\in C^{k+2,\alpha }(\Sigma )$ telle que
\begin{equation}
{\tilde L}[w]u=z,
\end{equation}
eu \'egard au fait que ${\tilde L}[w]$ est inversible. Or $z\in A^{k,\alpha }(\Sigma
)$, donc $G^{kl}D_{kl}z=0$. Ainsi, d\'erivons deux fois $(5.1)$ suivant des directions
horizontales $k$ et $l$ et saturons par $G^{kl}$, compte tenu de l'expression
$(2.15)$ donnant les composantes de la courbure de $D$ et puisque les d\'eriv\'ees covariantes horizontales de $w$ sont identiquement nulles, on obtient
$$G^{ij}D_{ij}(G^{kl}D_{kl}u)+B^{\alpha \beta }(w)D_{\alpha \beta }(G^{kl}D_{kl}u)-G^{kl}D_{kl}u=0$$
ce qui signifie que $G^{kl}D_{kl}u$ est dans le noyau de ${\tilde L}[w]$, donc, ce
terme est identiquement nul et par suite $(5.1)$ se r\'eduit \`a la forme suivante:
\begin{equation}
L[w](u)=z.
\end{equation}
D\'erivons l'\'equation $(5.2)$ par rapport \`a une direction horizontale $i\in \{1,...,n\}$, on obtient
$$B^{\alpha \beta }(w)D_{i\alpha \beta}u+D_{i}B^{\alpha \beta
}(w)D_{\alpha \beta }u-D_{i}u=0.$$
Saturons par $D^{i}u$, on voit que,
$$B^{\alpha \beta}D_{i\alpha \beta }uD^{i}u+D_{i}B^{\alpha \beta
}(w)D_{\alpha \beta}uD^{i}u-D_{i}uD^{i}u=0.$$
Tenons compte du fait que $w\in A^{k+1,\alpha }(\Sigma )$, qui implique que
$$D_{i}B^{\alpha \beta }(w)D_{\alpha \beta }uD^{i}u=0,$$
on aboutit \`a l'\'egalit\'e suivante
\begin{equation}
B^{\alpha \beta }D_{i\alpha \beta}uD^{i}u-{\tilde v}_{2}=0
\end{equation}
o\`u ${\tilde v}_2=D_{i}uD^{i}u$ d\'esigne le carr\'e de la norme de la composante
horizontale du gradient de $u$. D'apr\`es $(2.15)$, on a $D_{i\alpha \beta }u=D_{\alpha \beta i}u$ et donc l'\'egalit\'e $(5.3)$ s'\'ecrit
\begin{equation}
B^{\alpha \beta}D_{\alpha \beta i}uD^{i}u-{\tilde v}_2=0
\end{equation}
D'autre part, en un point $\xi _{0}$ o\`u ${\tilde v}_2$ atteint son maximum, on a: $B^{\alpha \beta }D_{\alpha \beta }{\tilde v}_2\leq 0$. Une in\'egalit\'e que se traduit par
$$B^{\alpha \beta }D_{\alpha \beta i}uD^{i}u+B^{\alpha \beta
}D_{\alpha i}uD_{\beta }^{i}u\leq 0$$
ce qui, compte tenu de $(5.4)$ et puisque le terme $B^{\alpha \beta }D_{\alpha  i}uD_{\beta }^{i}u$ est positif, implique que ${\tilde v}_{2}(\xi _{0})=0$. Ainsi ${\tilde v}_2$ est partout nul et le lemme est prouv\'e.

\vskip6mm

A pr\'esent, pour tout $t\in [0,1]$ et pour tout $w\in A^{1,\alpha }(\Sigma )$, on
pose $T_{t}w=u_{t}$, o\`u $u_{t}\in A^{2,\alpha }(\Sigma )$ est l'unique solution de
\begin{equation}
L[w]u_{t}=t\left[-w+(m-1)[1+v_1(w)]-[1+v_1(w)]^{\frac{3}{2}}e^{w}K(e^{w}\xi
)\right].
\end{equation}
Le lemme 3 montre que $u_{t}\in A^{2,\alpha }(\Sigma )$ et la th\'eorie elliptique
implique que
\begin{equation}
\Vert u_{t}\Vert _{\mathscr{C}^{2,\alpha }(\Sigma )}\leq Cste.
\end{equation}
La famille d'op\'erateurs $T_{t}$ v\'erifie les propri\'et\'es suivantes:
\begin{enumerate}
\item[(i)] Pour tout $w\in A^{1,\alpha }(\Sigma )$, $T_{0}w=0$.
\item[(ii)] A $t$ fix\'e, $w\mapsto T_{t}w$ est un op\'erateur compact de $A^{1,\alpha
}(\Sigma )$ dans $A^{1,\alpha }(\Sigma )$. En effet, soit $\Gamma \subset A^{1,\alpha
}(\Sigma )$ qu'on suppose born\'e dans $A^{1,\alpha }(\Sigma )$. Donc, d'apr\`es $(5.6)$, $\{T_{t}w\mid w\in \Gamma \}$ est un sous ensemble de $A^{2,\alpha }(\Sigma )$ qui
est uniform\'ement born\'e dans $C^{2,\alpha }(\Sigma )$. L'ensemble
$C^{2,\alpha }(\Sigma )$ \'etant un compact de $C^{1,\alpha }(\Sigma )$, on en
d\'eduit que $\{ T_{t}w\mid w\in \Gamma \}$ est un sous ensemble de $A^{1,\alpha
}(\Sigma )$ relativement compact.   
\item[(iii)] A $w$ fix\'e, $t\mapsto T_{t}w$ est continue, sinon, il existerait
une suite $(t_{i})$ dans l'intervalle $[0,1]$ et un r\'eel $\epsilon >0$ tels
que $\displaystyle \lim _{i\rightarrow \infty }t_{i}=t$ et, pour $i$ assez grand,
\begin{equation}
\Vert {\tilde u_{t_{i}}}-{\tilde u_{t}}\Vert _{C^{1,\alpha }(\Sigma )}>\epsilon .
\end{equation}
L'estim\'ee $(5.6)$ implique que la suite $({\tilde u_{t_{i}}})$ est born\'ee dans
$C^{2,\alpha }(\Sigma )$. L'injection de $C^{2,\alpha }(\Sigma )$ dans $C^{1,\alpha
}(\Sigma )$ \'etant compacte, donc quitte \`a en extraire une sous-suite, on peut supposer que $({\tilde u_{t_{i}}})$ converge dans $C^{1,\alpha }( \Sigma )$. D'apr\`es
$(5.5)$, et puisque ${\tilde u_{t_{i}}}$ est de classe $C^{2,\alpha }(\Sigma )$,
sa limite est ${\tilde u_{t}}$. Il en d\'ecoule que, pour \it i \rm assez grand,
$$\Vert {\tilde u_{t_{i}}}-{\tilde u_{t}}\Vert _{C^{1,\alpha }(\Sigma )}<\epsilon /2.$$
Ceci contredit $(5.7)$.
\end{enumerate}
 
\vskip2mm

Compte tenu des propri\'et\'es (i), (ii) et (iii), l'op\'erateur $T$ d\'efini en posant
$$T(t,w)=w-T_{t}w\mbox{ pour }(t,w)\in [0,1]\times A^{1,\alpha }(\Sigma )$$
est une d\'eformation compacte de l'identit\'e de $A^{1,\alpha }(\Sigma )$. Notre but est de r\'esoudre l'\'equation
\begin{equation}
T(t,u)=0,\ {\rm pour\ \it }t\in [0,1],
\end{equation} 
dans $A^{1,\alpha }(\Sigma )$. Le lemme 2 et la th\'eorie classique des \'equations
elliptiques, cf. [6], implique qu'une telle solution $u$ est dans $A^{2,\alpha }(\Sigma )$ et qu'il existe un r\'eel  $R>0$ telle que
\begin{equation}
\Vert u\Vert _{\mathscr{C}^{1,\alpha }(\Sigma )}<R.
\end{equation}
On note alors  $$B_{R}=\{ u\in A^{1,\alpha }(\Sigma )mid \Vert u\Vert _{\mathscr{C}^{1,\alpha }(\Sigma )}<R\}.$$
Compte tenu de $(5.9$), l'\'equation $(5.8)$ n'admet pas de solution $u$ sur le bord $\partial B_{R}$ de $B_{R}$. La d\'eformation $T$ est donc une homotopie compacte sur le
bord de $B_{R}$. En cons\'equence, d'apr\`es [3], (th\'eor\`eme 5.3.14), le
d\'egr\'e de Leray-Schauder de $T$ en $0$ relativement \`a $B_{R}$ ne d\'epend pas de $t$. Ainsi, pout tout $t\in [0,1]$, on a: 
$$d(T(t,.),0,B_{R})=d(T(0,.),0,B_{R})=d(Id,0,B_{R})=1.$$
En particulier, l'op\'erateur $T_{1}$ admet un point fixe $u\in A^{1,\alpha }(\Sigma
)$, solution de classe $A^{2,\alpha }(\Sigma )$ de l'\'equation 
$$\sum _{n+1\leq \alpha ,\beta \leq n+m-1}B^{\alpha \beta }D_{\alpha \beta
}u=(m-1)(1+v_{1})-(1+v_{1})^{\frac{3}{2}}e^{u}K(e^{u}\xi ).$$ 
Or $v_{2}(u)$ est partout nul, notre solution v\'erifie en fait l'\'equation
$$A^{ab}{\tilde D}_{ab}u=-v_{2}+(m-1)(1+v)-(1+v)^{\frac{3}{2}}e^{u}K(e^{u}\xi ).$$
Une \'equation qui est uniform\'ement elliptique ce qui permet de conclure en ce qui concerne la r\'egularit\'e $\mathscr{C}^{\infty }(\Sigma )$ de toute solution de classe $\mathscr{C}^{2,\alpha }(\Sigma )$.  

\vskip6mm

\subsection{Preuve du th\'eor\`eme 3}

\vskip4mm

La fonction $K\in \mathscr{C}^{\infty }(E_{*})$ est donn\'ee telle qu'il existe deux r\'eels $r_{1}$ et $r_{2}$ tels que $0<r_{1}\leq 1\leq r_{2}$ et
\begin{equation}
K(\xi )>\frac{m-1}{\Vert \xi \Vert }\ \mbox {si}\ \Vert \xi \Vert <r_{1},\
K(\xi )<\frac{m-1}{\Vert \xi \Vert }\ \mbox {si}\ \Vert \xi \Vert >r_{2}.
\end{equation}
On cherche \`a r\'esoudre dans $\displaystyle \Sigma $ l'\'equation elliptique suivante:
\begin{equation}
A^{ab}(u){\tilde D}_{ab}u=-v_{2}+(m-1)(1+v)-(1+v)^{\frac{3}{2}}e^{u}K(e^{u}\xi
).
\end{equation}

\vskip1mm

Tout d'abord, on note $\displaystyle \Sigma '=\{\xi \in E\mid r_{1}\leq \Vert \xi \Vert
\leq r_{2}\}$, $\nu $ le champ radial unitaire et $r$ la fonction $r(\xi
)=\Vert \xi \Vert $. Soit $\displaystyle w\in \mathscr{C}^{\infty }(\Sigma ')$, on d\'esigne par $w_{1}$ la restriction de $w$ \`a $\Sigma $ que l'on prolonge en une fonction radialement constante et pour $\displaystyle u\in \mathscr{C}^{\infty }(\Sigma ')$, on note
$\displaystyle {\tilde D}^{a}u=e^{\mu _{a}w_{1}}D^{a}u$, $\displaystyle
{\tilde D}^{ab}u=e^{(\mu _{a}+\mu _{b})w_{1}}D^{ab}u$,
$$\begin{array}{l}\displaystyle v_{1}(u)=(1-\mu _{a}){\tilde D}_{a}u{\tilde D}^{a}u,\
v_{2}(u)=\mu _{a}{\tilde D}_{a}u{\tilde D}^{a}u,\ {\tilde v}(u)=r^{2}v_{1}(u)+v_{2}(u),\\ \\ \displaystyle F(w)=-v_{2}(w_{1})+(m-1)\left[1+{\tilde
v}(w_{1})\right]-\left[1+{\tilde v}(w_{1})\right]^{3/2}e^{w_{1}}
K\left(e^{w_{1}}\frac{\xi }{\Vert \xi\Vert }\right),\\ \\
\displaystyle A^{ab}(w)=\left(1+{\tilde v}(w_{1})\right)G^{ab}-r^{2-\mu _{a}-\mu
_{b}}{\tilde D}^{a}w_{1}{\tilde D}^{b}w_{1},\end{array}$$
et $\displaystyle B(w)=(m-1)\left[1+{\tilde v}(w_{1})\right]-r^{2}v_{1}(w_{1})$.
Si $\displaystyle t\in [0,1]$ et $\displaystyle w\in \mathscr{C}^{\infty
}(\Sigma ')$, on d\'esigne par ${\tilde u_{t}}$ l'unique solution du probl\`eme de Neumann
\begin{equation}
\left\{\begin{array}{c}\displaystyle D_{\nu \nu }u+A^{ab}(w)r^{2-\mu _{a}-\mu
_{b}}{\tilde D}_{ab}u+\mu _{a}\log (r){\tilde D}^{a}_{a}u-u=t\Big[-w_{1}+F(w)\Big]\\ \\ \displaystyle +D_{\nu \nu }w+rB(w)D_{\nu }w-t\alpha \log (r)\sqrt{{\tilde v}(w_{1})}\ \mbox{ dans }\  \Sigma '\\ \\ \displaystyle D_{\nu }u=0\ \mbox{ sur }\ \ \partial \Sigma ',\end{array}
\right.
\end{equation} 
o\`u 
$$\alpha =2\sup _{\xi \in \Sigma }\sup _{\Sigma '}\Big\vert r\frac{\partial[rK(r\xi )]}{\partial r}\Big\vert .$$
La preuve de l'existence d'une solution ${\tilde u_{t}}\in C^{2,\alpha }(\Sigma ')$ de $(5.12)$ ainsi que de l'existence d'une constante positive $C$ telle que
$$\Vert {\tilde u_{t}}\Vert _{C^{2,\alpha }(\Sigma ')}\leq C$$
est tr\`es classique. On renvoie \`a [1] et [4], voir aussi [5], quant \`a l'unicit\'e, elle d\'ecoule du principe du maximum. L'\'equation $(5.12)$ \'etant elliptique et toutes les donn\'ees sont de classe $\mathscr{C}^{\infty }$, la solution ${\tilde u_{t}}$ est de classe $C^{\infty }(\Sigma ')$ par r\'egularit\'e. En particulier, \'etant donn\'ee une partie ${\cal B}$ born\'ee de $\mathscr{C}^{\infty }(\Sigma ')$, il existe ${\tilde u_{t}}\in \mathscr{C}^{\infty }(\Sigma ')$ ainsi qu'une suite
de r\'eels positifs $(C_{k})_{k\geq 0}$ telle que, quel que soit 
$(t,w)\in [0,1]\times {\cal B}$ et pour tout entier $k\geq 0$, on ait : 
\begin{equation}
\Vert {\tilde u_{t}}\Vert _{\mathscr{C}^{k}(\Sigma
')}\leq C_{k}.
\end{equation}
Il en d\'ecoule que l'op\'erateur $T_{t}$ d\'efini sur $\mathscr{C}^{\infty }(\Sigma ')$ par $T_{t}w={\tilde u_{t}}$ est compact. Il en est de m\^eme de l'op\'erateur
$T$ d\'efini de $[0,1]\times \mathscr{C}^{\infty }(\Sigma ')$ vers $\mathscr{C}^{\infty }(\Sigma ')$ par 
\begin{equation}
T(t,w)=w-T_{t}w.
\end{equation}
Notre objectif est de montrer que l'\'equation
\begin{equation}
T(t,u)=0,\mbox{ pour }t\in [0,1],
\end{equation}
admet une solution dans $\mathscr{C}^{\infty }(\Sigma ')$. Si une telle solution existe, celle-ci, not\'ee ${\tilde u_{t}}$, est une constante radiale. En effet, ${\tilde u_{t}}$ v\'erifie le syst\`eme suivant: 
\begin{equation}
\left\{\begin{array}{c}A^{ab}({\tilde u_{t}})r^{2-\mu _{a}-\mu
_{b}}{\tilde
D}_{ab}{\tilde u_{t}}+\mu _{a}\log (r){\tilde D}^{a}_{a}{\tilde
u_{t}}-{\tilde
u_{t}}=t\Big[-u_{t}+F({\tilde u_{t}})\Big]\\ \\
\displaystyle +rB({\tilde u_{t}})D_{\nu }{\tilde
u_{t}}-t\alpha \log (r)\sqrt{{\tilde v}(u_{t})}\ \mbox{ dans }\ \Sigma
'\\ \\ \displaystyle D_{\nu }{\tilde u_{t}}=0\ \mbox{ sur }\ \ \partial
\Sigma ',\end{array} \right.
\end{equation}
o\`u l'on a not\'e par $u_{t}$ le prolongement en une constante radiale de la restriction \`a $\Sigma $ de ${\tilde u_{t}}$.

\vskip3mm

D\'erivons radialement l'\'equation $(5.16)$ et multiplions par $r$ l'\'equation ainsi
obtenue. Du fait que $D_{\nu }r=1$ et puisque $u_{t}$ est une constante radiale, il en d\'ecoule que
\begin{equation}
\begin{array}{c}\displaystyle rD_{\nu }\Big[A^{ab}({\tilde
u_{t}})\Big]r^{2-\mu
_{a}-\mu _{b}}{\tilde D}_{ab}{\tilde u_{t}}+(2-\mu _{a}-\mu
_{b})r^{2-\mu
_{a}-\mu _{b}}A^{ab}({\tilde u_{t}}){\tilde D}_{ab}{\tilde
u_{t}}\\ \\ \displaystyle
+A^{ab}({\tilde u_{t}})r^{2-\mu _{a}-\mu
_{b}}r{\tilde D}_{\nu ab}{\tilde u_{t}}+\mu _{a}{\tilde
D}^{a}_{a}{\tilde
u_{t}}+\mu _{a}\log (r)G^{ab}r{\tilde D}_{\nu ab}{\tilde u_{t}}-rD_{\nu
}{\tilde u_{t}}\\ \\ \displaystyle =trD_{\nu }F({\tilde u_{t}})+rB({\tilde
u_{t}})D_{\nu }{\tilde
u_{t}}+rD_{\nu }{\tilde u_{t}}rD_{\nu }\Big[B({\tilde
u_{t}})\Big]+B({\tilde
u_{t}})r^{2}D_{\nu \nu }{\tilde u_{t}} \\ \\
\displaystyle -t\alpha \sqrt{{\tilde
v}(u_{t})}-t\alpha \log (r)rD_{\nu
}\sqrt{r^{2}v_{1}(u_{t})+v_{2}(u_{t})}.\end{array}
\end{equation}
Remarquons d'abord que par un calcul analogue \`a celui de la seconde section, on montre que
\begin{equation}
D_{e_{a}}(r\nu )=(1-\mu _{a})e_{a},
\end{equation}
o\`u $\mu _{a}$ vaut $1$ ou $0$ selon que la direction $e_{a}$ est horizontale ou verticale. Ainsi, usons de la d\'efinition de la d\'eriv\'ee covariante, on peut \'ecrire: 
$$D_{e_{a}}(rD_{\nu }u)=D^{2}u(e_{a},r\nu )+Du\Big(D_{e_{a}}(r\nu )\Big)=rD_{\nu
}(D_{a}u)+(1-\mu _{a})D_{a}u.$$
On en d\'eduit que
$$rD_{\nu}(D_{a}u)=D_{e_{a}}(rD_{\nu }u)-(1-\mu _{a})D_{a}u$$
et, en particulier, pour une constante radiale $u$, on a :
\begin{equation}
\begin{array}{c}\displaystyle rD_{\nu }v_{1}(u)=2(1-\mu
_{a})D_{a}(rD_{\nu }u) D^{a}u-2v_{1}(u)=-2v_{1}(u),\\ \\ \displaystyle
rD_{\nu }v_{2}(u)=2\mu _{a}{\tilde D}_{a}(rD_{\nu }u) {\tilde
D}^{a}u+2v_{2}(u)rD_{\nu }u=0\end{array}
\end{equation}
Combinons cette relation avec le fait que $u_{t}$ est une constante radiale, on
obtient:
$$rD_{\nu }\left[r^{2}v_{1}(u_{t})+v_{2}(u_{t})\right]=0.$$
Ainsi, l'\'equation $(5.17)$ se r\'eduit \`a la suivante :
\begin{equation}
\begin{array}{c}\displaystyle
(2-\mu _{a}-\mu _{b})r^{2-\mu _{a}-\mu _{b}}A^{ab}({\tilde
u_{t}}){\tilde D}_{ab}{\tilde u_{t}} +A^{ab}({\tilde u_{t}})r^{2-\mu _{a}-\mu
_{b}}r{\tilde D}_{\nu ab}{\tilde u_{t}}\\ \\ \displaystyle +\mu
_{a}{\tilde D}^{a}_{a}{\tilde u_{t}}+\mu _{a}\log (r)G^{ab}r{\tilde D}_{\nu
ab}{\tilde u_{t}}-rD_{\nu }{\tilde u_{t}}=rB({\tilde u_{t}})D_{\nu }{\tilde
u_{t}}\\ \\ \displaystyle +B({\tilde u_{t}})r^{2}D_{\nu \nu
}{\tilde u_{t}}-t\alpha \sqrt{r^{2}v_{1}(u_{t})+v_{2}(u_{t})}.\end{array}
\end{equation}
D'autre part, la d\'efinition de la d\'eriv\'ee covariante nous permet d'\'ecrire que  $$D_{ab}(rD_{\nu}u)=D^{2}(rD_{\nu }u)(e_{a},e_{b})=e_{a}\Big[D(rD_{\nu
}u)(e_{b})\Big]-D(rD_{\nu }u)(D_{e_{a}}e_{b}).$$
D'o\`u
$$\begin{array}{c}\displaystyle D_{ab}(rD_{\nu }u)=e_{a}\left[D_{e_{b}}\left(Du(r\nu
)\right)\right]-\left(D_{e_{a}}e_{b}\right)Du(r\nu )\\ \\\displaystyle
=D_{e_{a}}(D^{2}u(e_{b},r\nu )+D_{e_{a}}Du\left(D_{e_{b}}(r\nu
)\right)\\ \\ \displaystyle -D^{2}u\left(D_{e_{a}}e_{b},r\nu
\right)-Du\left(D_{D_{e_{a}}e_{b}}(r\nu )\right)\end{array}$$
et par suite  
$$\begin{array}{c}\displaystyle D_{ab}(rD_{\nu
}u)=D^{3}u(e_{a},e_{b},r\nu
)+D^{2}u\Big(e_{b},D_{e_{a}}(r\nu )\Big)
+D^{2}u\Big(e_{a},D_{e_{b}}(r\nu )\Big)\\ \\ \displaystyle
+Du\Big(D_{e_{a}}(D_{e_{b}}(r\nu
))\Big)-Du\Big(D_{D_{e_{a}}e_{b}}(r\nu )\Big).\end{array}$$
Or, la relation $(5.18)$ permet d'\'ecrire que
$$Du\left(D_{e_{a}}(D_{e_{b}}(r\nu ))\right)=(1-\mu
_{b})Du\left(D_{e_{a}}e_{b}\right)$$
et 
$$Du\left(D_{D_{e_{a}}e_{b}}(r\nu )\right)=(1-\mu
_{b})Du\left(D_{e_{a}}e_{b}\right).$$
Ainsi, usons \`a nouveau de $(5.18)$, on obtient : 
$$D_{ab}(rD_{\nu }u)=D^{3}u(e_{a},e_{b},r\nu )+(2-\mu _{a}-\mu _{b})D^{2}u(e_{a},e_{b})$$
et d'apr\`es la d\'efinition de la connexion $D$ et celle de la d\'eriv\'ee covariante, on voit que
\begin{equation}
D_{ab}(rD_{\nu }u)=rD_{\nu }(D_{ab}u)+(2-\mu
_{a}-\mu _{b})D_{ab}u.
\end{equation}
Eu \'egard \`a la relation $(2.6)$ de la secande section donnant l'expression de la courbure de $D$ et dont d\'ecoule l'\'egalit\'e $\displaystyle D_{ab\nu }u=D_{\nu ab}u$. Reportons $(5.21)$ dans $(5.20)$, on obtient :
\begin{equation}
\begin{array}{c}\displaystyle A^{ab}({\tilde u_{t}})r^{2-\mu _{a}-\mu
_{b}}{\tilde
D}_{ab}(rD_{\nu }{\tilde u_{t}})+\mu _{a}{\tilde
D}^{a}_{a}{\tilde u_{t}}+\mu
_{a}\log (r){\tilde D}^{a}_{a}(rD_{\nu }{\tilde
u_{t}})\\ \\ \displaystyle -rD_{\nu }{\tilde
u_{t}}=rB({\tilde
u_{t}})D_{\nu }{\tilde u_{t}}+B({\tilde u_{t}})r^{2}D_{\nu \nu
}{\tilde
u_{t}}-t\alpha \sqrt{{\tilde v}(u_{t})}.\end{array}
\end{equation}

\vskip0mm

Rappellons que l'\'equation de Gauss et les calculs de la seconde section impliquent que pour deux fonctions $u_{1},u_{2}\in \mathscr{C}^{2}(\Sigma ')$ ayant les m\^emes valeurs sur $\Sigma $, on a :
\begin{equation}
D_{ab}u_{1}=D_{ab}u_{2}+(1-\mu _{a})G_{ab}D_{\nu}(u_{1}-u_{2}),\ a,b\leq
n+m-1.
\end{equation}
Or $u_{t}=({\tilde u_{t}})_{\vert \Sigma }$. Donc, tenons compte de $(5.23)$ et du fait que $D_{\nu }{\tilde u_{t}}$ est nul sur $\Sigma $, $(5.22)$ implique la relation suivante : 
\begin{equation}
\mu _{a}{\tilde D}^{a}_{a}u_{t}=-t\alpha 
\sqrt{v_{1}(u_{t})+v_{2}(u_{t})},\mbox{ partout dans }\Sigma .
\end{equation}
D'autre part, restreignons $(5.16)$ \`a $\Sigma $, l'usage de $(5.23)$ montre que la fonction $u_{t}$ v\'erifie l'\'equation suivante :
\begin{equation}
A^{ab}(u_{t}){\tilde D}_{ab}u_{t}-u_{t}=t\Big[-u_{t}+F(u_{t})\Big],\mbox{ partout
dans }\Sigma .
\end{equation}
La fonction $u_{t}$ \'etant une constante radiale, une combinaison des \'equations $(5.24)$ et $(5.25)$, montre que $u_{t}$ est une autre solution de $(5.16)$. En
cons\'equence, on voit que
$$\left\{ \begin{array}{c}\displaystyle A^{ab}({\tilde u_{t}})r^{2-\mu _{a}-\mu
_{b}}{\tilde
D}_{ab}({\tilde u_{t}}-u_{t})+\mu _{a}\log (r){\tilde
D}^{a}_{a}({\tilde u_{t}}-u_{t})\\ \\
\displaystyle -({\tilde
u_{t}}-u_{t})=rB({\tilde u_{t}})D_{\nu }({\tilde u_{t}}-u_{t})\ {\rm
dans\it
}\ \Sigma '\\ \\ \displaystyle D_{\nu }({\tilde u_{t}}-u_{t})=0\ \
{\rm sur\it }\ \ \partial
\Sigma '.\end{array}\right. $$
Le principe du maximum montre que ${\tilde u_{t}}=u_{t}$ partout dans $\Sigma
'$. Ceci montre que ${\tilde u_{t}}$ est une constante radiale.

\vskip3mm

A pr\'esent, on montre que $u_{t}$ est estim\'ee \`a priori $\mathscr{C}^{0}(\Sigma )$. Une majoration a priori se d\'eduit imm\'ediatement du principe du maximum. Soit $\xi \in \Sigma $ un point o\`u $u_{t}$ atteint son maximum. Si $u_{t}(\xi )>\log (r_{2})$, l'hypoth\`ese de croissance $(5.10)$ faite sur $K$ combin\'ee avec l'\'equation $(5.25)$ implique qu'au point $\xi $, on aura 
$$0\geq {\tilde D}_{a}^{a}u_{t}=(1-t)u_{t}+t\left[ (m-1)-e^{u_{t}}K(e^{u_{t}}\xi
)\right]>0$$
ce qui constitue une contradiction, eu \'egard au fait que $\log (r_{2})\geq 0$. La minoration $u_{t}\geq \log (r_{1})$ s'obtient par analogie en consid\'erant un point o\`u $u_{t}$ atteint son minimum.

\vskip3mm

La fonction $u_{t}$ v\'erifiant les \'equations $(5.24)$ et $(5.25)$, d'apr\`es le lemme 1 elle est estim\'ee a priori dans $\mathscr{C}^{1}(\Sigma )$. La th\'eorie classique des \'equations uniform\'ement elliptiques, voir \`a titre d'exemple [1], [4], [6], montre comment obtenir l'estim\'ee $\mathscr{C}^{1,\alpha }$. Plus
pr\'ecis\'ement, on peut, par simple modification de la preuve du th\'eor\`eme 13.6 page 328 de Gilbarg et Tr\"udinger [4], conclure \`a l'existence d'un
r\'eel $c_{1}>0$ tel que 
$$\Vert u_{t}\Vert _{\mathscr{C}^{1,\alpha }(\Sigma )}<c_{1}.$$
La fonction $u_{t}$ \'etant une solution de l'\'equation $(5.25)$, les in\'egalit\'es de Schauder impliquent l'existence d'une constante positive $C$ telle que
$$\Vert u_{t}\Vert _{\mathscr{C}^{2,\alpha }(\Sigma )}<C\left[\Vert
u_{t}\Vert _{\mathscr{C}^{0}(\Sigma )}+\Vert F\Vert _{\mathscr{C}^{0,\alpha }(\Sigma
')}\right]$$
et donc il existe une constante positive $c_{2}$ telle que
$$\Vert u_{t}\Vert _{\mathscr{C}^{2,\alpha }(\Sigma
)}<c_{2}.$$
Supposons donc que que pour tout $s\leq k$, pour un $k\geq 2$, on ait:
$$\Vert u_{t}\Vert
_{\mathscr{C}^{s,\alpha }(\Sigma )}<c_{s}.$$
L'\'equation obtenue en d\'erivant $(k-1)$ fois l'\'equation $(5.25)$ s'\'ecrit localement sous la forme
$$A^{ab}{\tilde
D}_{ab}(D_{i_{1}i_{2}...i_{k-1}}u_{t})=H_{i_{1}i_{2}...i_{k-1}},$$
o\`u le second membre $\displaystyle H_{i_{1}i_{2}...i_{k-1}}$ ne d\'epend que des
d\'eriv\'ees covariantes de $u_{t}$ d'ordre $\leq k$. Il est donc born\'e dans
$\mathscr{C}^{0,\alpha }(\Sigma ')$ et il en d\'ecoule, d'apr\`es les in\'egalit\'es de Schauder, que $\displaystyle \Vert D^{(k-1)}u_{t}\Vert _{\mathscr{C}^{2,\alpha }(\Sigma )}<c_{k+1}$, et par suite
$$\Vert u_{t}\Vert _{\mathscr{C}^{k+1,\alpha
}(\Sigma )}<c_{k+1}.$$
On vient ainsi d'\'etablir par r\'ecurrence l'existence de r\'eels positifs $a_{k}$ tels que
\begin{equation}
\Vert
u_{t}\Vert _{\mathscr{C}^{k}(\Sigma ')}<a_{k}.
\end{equation}
Eu \'egard au fait que $u_{t}$ est radialement constante. On note $B$ la pseudo-boule d\'efinie par 
$$B=\{ u\in \mathscr{C}^{\infty
}(\Sigma ')\mid \Vert u\Vert _{\mathscr{C}^{k}(\Sigma ')}<a_{k}\}.$$
Compte tenu de $(5.26)$, l'\'equation $(5.15)$ n'admet pas de solution sur le bord de $B$. La d\'eformation $T$ est donc une homotopie compacte sur le bord de $B_{R}$. En cons\'equence, d'apr\`es le th\'eor\`eme de Nagumo [7], le d\'egr\'e de $T$ en $0$
relativement \`a $B$ ne d\'epend pas de $t$. Ainsi, pout tout $t\in [0,1]$, on : 
\begin{equation}
d(T(t,.),0,B)=d(T(0,.),0,B)=\gamma .
\end{equation}
Or, pour $t=0$, la fonction $u_{0}=0$ est l'unique solution de $(5.15)$ et l'on montre
ais\'ement que pour $\displaystyle w\in \mathscr{C}^{\infty }(\Sigma ')$,$(d_{u_{0}}T_{0})(w)=u$, o\`u $u$ est l'unique solution du probl\`eme suivant :
\begin{equation}
\left\{ \begin{array}{c}\displaystyle
D_{\nu \nu }u+r^{2-2\mu _{a}}D^{a}_{a}u+\mu
_{a}\log (r)D^{a}_{a}u-u=D_{\nu
\nu }w\\ \\ \displaystyle +(m-1)rD_{\nu }w\ {\rm dans\it }\ \Sigma
'\\ \\ \displaystyle D_{\nu }u=0\ \ {\rm sur\it }\ \ \partial \Sigma '.\end{array}
\right.
\end{equation}
A pr\'esent, si $w\in \ker [\mbox{id}-(d_{u_{0}}T_{0})]$, celle-ci sera solution du probl\`eme suivant :
$$\left\{ \begin{array}{c}\displaystyle r^{2-2\mu
_{a}}D^{a}_{a}w+\mu
_{a}\log (r)D^{a}_{a}w-w=(m-1)rD_{\nu }w\ {\rm dans\it
}\ \Sigma '\\ \\ \displaystyle D_{\nu
}w=0\ \ {\rm sur\it }\ \ \partial
\Sigma '.\end{array} \right. $$
Un raisonnement analogue \`a celui qui pr\'ec\`ede montre que $w$ est une constante radiale. Celle-ci \'etant identiquement nulle sur $\Sigma $, donc partout nulle et par suite
\begin{equation}
\ker [\mbox{id}-(d_{u_{0}}T_{0})]=\{0\}.
\end{equation}
Or, d'apr\`es sa d\'efinition et en raisonnant comme auparavant, on v\'erifie que $d_{u_{0}}T_{0}$ est un op\'erateur compact de $\mathscr{C}^{\infty }(\Sigma ')$ vers $B^{1,\beta }(\Sigma ')$, $0<\beta <\alpha $, et donc, tenons compte de $(5.29)$, le
th\'eor\`eme d'alternative de Fredholm implique que l'op\'erateur $\displaystyle
\mbox{id}-(d_{u_{0}}T_{0})$ est inversible. En cons\'equence $0$ est un point r\'egulier pour $\mbox{id}-T_{0}$. Ainsi, dans $(5.27)$, $\gamma =\pm 1$ et en particulier, l'op\'erateur $T_{1}$ admet un point fixe de classe $\mathscr{C}^{\infty }(\Sigma ')$. Celle-ci est une constante radiale qui est une solution de $(5.11)$ dans $\Sigma $. Ceci ach\`eve la preuve du th\'eor\`eme. 

\vskip6mm

\subsection{Remarques}

\vskip4mm

\noindent\textbf{1-} Expliquons ici en quel sens l'hypoth\`ese de croissance des th\'eor\`emes 2 et 3 est la meilleure possible. Soit  $K\in \mathscr{C}^{\infty }(E_{*})$ une fonction strictement positive. On suppose qu'il existe un r\'eel $a\in ]0,1[$ tel que 
$$K(\xi )\leq \frac{a(m-1)}{\Vert \xi \Vert },\ \xi \in E_{*}.$$
Alors il n'existe pas de solution de classe $\mathscr{C}^{2}(\Sigma )$ de l'\'equation $(5.11)$ ci-dessus. En effet, si une telle solution existe, on aura
\begin{equation}
A^{ab}{\tilde D}_{ab}u\geq -v_{2}+(m-1)(1+v)-a(m-1)(1+v)^{\frac{3}{2}}.
\end{equation}
En un point $\xi \in \Sigma $ o\`u $u$ atteint son maximum, on a : $\displaystyle v(\xi )=0$ et dans un rep\`ere $G$-orthonorm\'e diagonalisant $(D_{ab}u(\xi ))$, $(5.30)$ s'\'ecrit au point $\xi $ sous la forme
\begin{equation}
\sum _{1\leq a\leq n+m-1}{\tilde D}_{aa}u\geq (m-1)(1-a).
\end{equation}
Or $D_{aa}u(\xi )\leq 0$. Reportons dans $(5.31)$, on voit que $1\leq a$ ce qui est contradictoire.

\vskip3mm

On obtient la m\^eme conclusion s'il existe un r\'eel $b>1$ tel que
$$K(\xi )\geq \frac{b(m-1)}{\Vert \xi \Vert },\ \xi \in E_{*}.$$

\vskip4mm

\noindent\textbf{2-} Consid\'erons une fonction $K$ v\'erifiant les hypoth\`eses du th\'eor\`eme 2 et supposons que pour tout $r\in [r_{1},r_{2}]$, il existe $\xi \in \Sigma $ telle que
\begin{equation}
K(r\xi )\neq (m-1)r^{-1}.
\end{equation}
D'apr\`es la preuve du th\'eor\`eme 2, il existe une solution $u_{1}\in \mathscr{C}^{\infty }(\Sigma )$ de l'\'equation $(5.11)$ ci-dessus \`a gradient horizontal identiquement nul, $v_{2}(u_{1})=0$, et telle que $\displaystyle r_{1}\leq e^{u_{1}}\leq r_{2}$. D'autre part, la preuve du th\'eor\`eme 3 donne une solution $u_{2}$ de cette m\^eme \'equation telle que $\displaystyle r_{1}\leq e^{u_{2}}\leq r_{2}$ et de sorte que pour un r\'eel $\alpha >0$, assez grand, on ait :
\begin{equation}
\mu _{a}e^{2u_{2}}D^{a}_{a}u_{2}=-\alpha \sqrt{v(u_{2})}.
\end{equation}
Les deux graphes radiaux donn\'es par ces deux solutions ne sont pas homoth\'etiques. Sinon, la fonction $u_{1}-u_{2}=Cste$, et donc $v_{2}(u_{1})=v_{2}(u_{2})=0$. Reportons dans $(5.33)$, et puisqu'alors $\displaystyle \mu _{a}e^{2u_{2}}D^{a}_{a}u_{2}$ est identiquement nul, il en d\'ecoule que $v_{1}(u_{2})=0$, et par suite $u_{2}$ est une constante. Il existe alors $r\in [r_{1},r_{2}]$ telle que $\displaystyle
e^{u_{2}}=r$. L'\'equation $(5.11)$ satisfaite par $u_{2}$ implique que $K(r\xi )=(m-1)r^{-1}$, quel que soit $\xi \in \Sigma $, contredisant ainsi $(5.32)$. 

\vskip8mm



\begin{thebibliography}{99}
%
\bibitem{label1}
S. Agmon: Lectures on Elliptic boundary value problems; Van
Nostrand,
Princeton, NJ 1965.
%
\bibitem{label2} I. Bakelman and B.
Kantor: Existence of spherically homeomorphic hypersurfaces
in euclidean
space with prescribed mean curvature; Geometry and Topology, Leningrad, 1,
1974,
3-10. 
%
\bibitem{label3} M. S. Berger: Nonlinearity and
functional analysis; Academic Press, New
York, San Francisco, London 1977,
Pure and Applied Mathematics 74. 
%
\bibitem{label4} D. Gilbarg and N.
Tr\"udinger: Elliptic partial differential equations of second
order; 2nd
ed., Springer, New York, 1977.   
%
\bibitem{label5} O.A. Ladyzhenskaya
and N.N. Ural'tseva: Linear and quasilinear elliptic
equations; Academic
press, New York and London, 1968.
%
\bibitem{label6} C. B. Morrey:
Multiple integrals in the calculus of variations; 130,
Springer-verlag,
Berlin, Heidelberg, New York, 1966.
%
\bibitem{label7} M. Nagumo: Degree
of mapping in convex linear topological spaces; Amer. J. Math.,
73, 1951, p.
497-511.
%
\bibitem{label8} A. E. Treibergs and S. W. Wei: Embedded
hypersurface with prescribed mean
curvature; J. Diff. Geo., 18, 1983,
513-521. 
%
\bibitem{label9} K. Yano and S. Ishihara: Tangent and
cotangent bundles; Marcel Dekker, New-York,
1973. 
%
\bibitem{label10}
K. Yano and S. Ishihara: Horizontal lifts of tensor fields and connections
to
tangent bundles; J. Math. and Mech., 16, 1967,
1015-1030.
\end{thebibliography}
\end{document}